\documentclass[12pt,reqno]{amsart}
\usepackage{latexsym,stmaryrd,amsmath,amsthm,amssymb,amsfonts}
\usepackage{setspace}
\usepackage{color}
\usepackage{tikz}
\usepackage{pgf,tikz}
\usepackage{mathrsfs}
\usetikzlibrary{arrows}
\numberwithin{equation}{section}
\newtheorem{theorem}{Theorem}[section]
\newtheorem{lem}[theorem]{Lemma}

\newtheorem{pro}[theorem]{Proposition}
\newtheorem{cor}[theorem]{Corollary}
\newtheorem{defi}[theorem]{Definition}
\newtheorem{rem}[theorem]{Remark}

\def\S2{\mathbb{S}^2}
\def\s{\,\,\,\,}

\def\endproof{$\hfill\Box$\\}
\def\R{\mathbb{R}}

\title{\bf Compactness of isospectral conformal metrics on 4-manifolds}

\author{Ke Xu}
\address{Department of Mathematical science\\Tsinghua University\\P. R. China}
\email{xuke16@mails.tsinghua.edu.cn}
\date{}
\begin{document}
\maketitle

\begin{abstract}
Let a sequence of conformal Riemannian metrics $\{g_k=u_k^2g_0\}$ be isospectral to $g_0$ over a compact boundaryless smooth 4-dimension manifold $(M,g_0)$. We prove that the subsequence of conformal factors $\{u_k\}$ converges to $u$ weakly in $W^{2,p}_{loc}(M\setminus \mathcal{S})$ for some $p<2$, where $\mathcal{S}$ is a finite set of points and $u\in W^{2,p}(M,g_0)$. Moreover, if the isospectral invariant $\frac{\int_M R(g_k)dV_{g_k}}{6\sqrt{\mathrm{Vol}(M,g_k)}}$ is strictly smaller than the Yamabe constant of the standard sphere $\mathbb{S}^4$, then the subsequence of distance functions $\{d_k\}$ defined by $\{g_k\}$ uniformly converges to $d_u$ and the subsequence of metric spaces $\{(M,d_k)\}$ converges to the metric space $(M,d_u)$ in the Gromov-Hausdorff topology, where $d_u$ is the distance function defined by $u^2g_0$.  
\end{abstract}

\section{Introduction}
Let $M$ be a compact $4$-manifold without boundary. A smooth Riemannian metric $g_0$ on $M$ determines a class of conformally equivalent metrics of the form $g=u^2g_0$, where $u$ is a smooth positive function.  Our main theorem is about the Gromov-Hausdorff compactness for metric spaces $(M,d)$ defined by $(M,g)$.
\begin{theorem}\label{main}
Let $(M,g_0)$ be a compact smooth $4$-dimensional Riemannian manifold without boundary. Assume $\{g_k=u_k^2g_0\}$ is a sequence of conformal metrics satisfying  
\begin{align}
\mathrm{Vol}(M,g_k)=1  
\end{align}
\begin{align} \label{1.1}
\int_M R(g_k)^{2}dV_{g_k}<A
\end{align}
for some positive constant $A$, where $R(g_k)$ is the scalar curvature of $g_k$.
 
Then $\{u_k\}$ weakly converges to $u$ in $W^{2,p}_{loc}(M\setminus \mathcal{S})$ for some $p<2$, where $\mathcal{S}$ is a finite set and $u\in W^{2,p}(M)$. Moreover, if

(1) The first eigenvalue of $\Delta_{g_k}$: $\lambda_1(\Delta_{g_k})\geq \Lambda$, for some positive constant $\Lambda$;

(2) $\liminf_{k\to \infty}\frac{1}{6}\int_M R(g_k)dV_{g_k}< Y(\mathbb{S}^4)$, $Y(\mathbb{S}^4)$ is the Yamabe constant of the standard sphere $\mathbb{S}^4$;
 
Then, after passing to a subsequence, a sequence of distance functions $\{d_k\}$ defined by $\{g_k\}$ uniformly converges to the distance function $d_u$ defined by $g_u=u^2g_0$. In other words, the sequence of metric spaces $\{(M,d_k)\}$ converges to the metric space $(M,d_u)$ in the Gromov-Hausdorff distance.
\end{theorem} 
\begin{rem}
In general, this theorem remains valid in any dimension $n\geq 3$. If we replace (\ref{1.1}) with
\begin{align}
\int_M |R(g_k)|^{\frac{n}{2}}dV_{g_k}<A
\end{align}
then we get $\{u_k\}$ weakly converges to $u$ in $W^{2,p}_{loc}(M\setminus \mathcal{S})$ and $u\in W^{2,p}(M)$, for $p<\frac{n}{2}$ . Moreover, if the constraint (2) is replaced by 

(2') $\liminf_{k\to \infty}\frac{n-2}{4(n-1)}\int_M R(g_k)dV_{g_k}< Y(\mathbb{S}^n)$, $Y(\mathbb{S}^n)$ is the Yamabe constant of the standard sphere $\mathbb{S}^n$.

Then  a sequence of distance functions $\{d_k\}$ defined by $\{g_k\}$ uniformly converges to the distance function $d_u$ defined by $g_u=u^\frac{4}{n-2}g_0$.
\end{rem}

The original motivation of the above result is the application to  isospectral conformal metrics. At first, we need to establish some basic conformal notions and recall some essential results on the spectral invariants. Let $(M,g_0)$ be a compact $n$-dim Riemannian manifold without boundary. Denote the scalar curvature of $g$ by $R(g)$ (or $R_g$). If we consider the metrics $g=u^\frac{4}{n-2}g_0$ for some positive smooth functions $u$, $u$ satisfy the following equation
$$
-\frac{4(n-1)}{n-2}\Delta u+R(g_0)u=R(g)u^\frac{n+2}{n-2}.
$$
Two Riemannian metrics $g$ and $g'$ on a compact manifold $M$ are said to be isospectral if their associated Laplacian operators have the identical spectrum, i.e., $\mathrm{Spec}(\Delta_g)=\mathrm{Spec}(\Delta_{g'})$. It is well-known that the heat kernel $H_t(x,y)$ has an eigenfunction expansion:
$$
H_t(x,y)=\sum_i e^{-\lambda_it}v_i(x)v_i(y),
$$
where $\{v_i(x)\}$ is an orthonormal basis for the eigenfunction space of $\Delta$. 
The trace of the heat kernel $e^{-t\Delta}$ has the known expansion as $t\to 0^+$
$$
\mathrm{Tr}(H_t)=\sum_i e^{-t\lambda_i}\sim (4\pi t)^{-\frac{n}{2}}(a_0+a_1t+a_2t^2+a_3t^3+\cdots),
$$
where each $a_i$ is the spectral invariant (cf.\cite{Gilkey}). The first several heat invariants $a_k$ are given by
\begin{align*}
a_0&=\int_MdV_g=\mathrm{Vol}(M,g);
\\a_1&=\frac{1}{6}\int_M R_g dV_g;
\\a_2&=\frac{1}{360}\int_M(5R_g^2-2|Ric_g|^2+2|Riem_g|^2)dV_g.
\end{align*}
And $R_g$, $Ric_g$ and $Riem_g$ denote the scalar curvature, Ricci curvature tensor and full curvature tensor respectively. Since we are going to work with conformal classes in dimension $4$, it's convenient for us to rewrite  $a_2$ as follows (cf. \cite{Xu1}).
$$
a_2=\frac{1}{180}\int_M(|W_g|^2+|B_g|^2+\frac{29}{12}R_g^2)dV_g,
$$
where $W_g$ and $B_g$ are the Weyl curvature tensor and the traceless Ricci curvature tensor. Note that $a_2$ implies that scalar curvature of isospectral conformal metrics are $L^2$-invariant in the $4$-dimensional case (see \cite{B},\cite{Xu1}). If $g_1$ is isospectral to $g_2$ and both of them are in the same conformal class, we have 
$$a_2(g_1)=a_2(g_2)\s and \s \int_{M}|W_{g_1}|^2dV_{g_1}=\int_M|W_{g_2}|^2dV_{g_2}.$$
The Gauss-Bonnet formula for four dimensional closed manifolds tells us that
$$
\int_M (R_{g_1}^2-12|B_{g_1}|^2)dV_{g_1}=\int_M (R_{g_2}^2-12|B_{g_2}|^2)dV_{g_2}.
$$ 
It implies that 
$$
\int_M R_{g_1}^2dV_{g_1}=\int_M R_{g_2}^2dV_{g_2}.
$$
Then, we have the following natural constraints for the isospectral conformal metrics $g=u^2g_0$: 
\begin{align*}
\int_M u^4dV_0&=a_0;
\\ \frac{1}{6}\int_M R_g u^4dV_0 &=a_1;
\\ \int_MR_g^2u^4dV_0&=a_2'\leq a_2.
\end{align*}

Now, we apply our main theorem to the isospectral conformal metrics and get
\begin{cor}
Assume a set of metric spaces $\{(M,d_k)\}$ is induced by an isospectral set of conformal Riemannian metrics $\{g_k=u_k^2g_0\}$ on a compact smooth $4$-manifold $(M,g_0)$ without boundary. Then $\{u_k\}$ converges to $u$ weakly in $W^{2,p}_{loc}(M\setminus \mathcal{S})$ and $u\in W^{2,p}(M)$ for some $p<2$, where $\mathcal{S}$ is a finite set. Moreover, if 
$$\frac{a_1}{\sqrt{a_0}}< Y(\mathbb{S}^4),$$
where $Y(\mathbb{S}^4)$ is the Yamabe constant of the standard sphere $\mathbb{S}^4$ and $a_0, a_1$ are the leading coefficients of the heat trace expansion $\mathrm{Tr}(H_t)$ at $t\to 0$. Then $\{d_k\}$ uniformly converges to the distance function $d_u$, where $d_u$ is defined by $g_u=u^2g_0$. In other words, $\{(M,d_k)\}$ converges to $(M,d_u)$ in the Gromov-Hausdorff topology.
\end{cor}

The compactness of isospectral metrics has been studied for a long time. For the case of compact surfaces, Osgood, Phillips and Sarnak\cite{OPS} showed that a set of isospectral metrics on a compact Riemann surface without boundary form a compact family in the $C^\infty$ topology. In the higher dimensions, if the isospectral metrics are restricted to the same conformal class, many good results can be derived from it. For three dimensional compact manifolds without boundary, A.Chang and P.Yang \cite{Chang-Yang2}proved that an isospectral set of conformal metrics is compact in the $C^\infty$-topology. As for the dimension greater than 3, Gursky \cite{Gursky} did some $C^\alpha$ compactness for manifolds with $L^p$ norm of the full curvature tensor bounds when $p>\frac{n}{2}$ and $0<\alpha<(2p-n)/p$. In particular, the $C^\infty$ topology compactness of isospectral conformal metrics over a closed $4$-manifold was presented under some extra constrains(see e.g.,\cite{L-W}, \cite{Xu1} ,\cite{Xu2}). While their results in the dimension $4$ are stronger than ours, we make only minimal assumptions on the isospectral invariant.

In the remainder of the introduction, we give an outline of the argument for our main theorem. We first show  $\varepsilon$-regularity of $\{u_k\}$ and define a finite set of points:
$$
\mathcal{S}=\{x\in M: \lim_{r\to 0}\liminf_{k\to \infty}\int_{B_r(x)} R_k^2dV_{g_k}>\varepsilon_0\},
$$
such that $\{u_k\}$ converges to $u$ weakly in $W^{2,p}_{loc}(M\setminus S)$, for some $p<2$. Then the removable singularities result can be proved by using the Three Circles Lemma in the end of section 2.

Then we focus on a sequence of distance functions $\{d_k\}$ defined by $\{g_k\}$ in section 3. According to our previous paper \cite{Dong-Li-Xu}, it asserts that there exists a finite set $\mathcal{S}'$, such that  $d_k$ converges to $d_u$ in $C^0_{loc}(M\setminus \mathcal{S}')$. If in addition, $u\neq 0$, then $d_u$ is a distance function defined by $g_u=u^2g_0$. Since singularities of $u$ are removable, $d_u$ is well-defined over $M\times M$. Moreover we show if $\lim_{r\to 0} \mathrm{diam}(B_r(P),g_k)+\mathrm{Vol}(B_r(P),g_k)=0$ for any $P\in \mathcal{S}'$, then $\{(M,d_k)\}$ converges to $(M,d_u)$ in the sense of Gromov-Hausdorff topology. 

In section 4, we firstly recall some basic terminology about bubble trees, one can refer to \cite{Chen-Li} and \cite{Parker}. Then we work on the neck analysis and show that the volume and the diameter of the neck region will vanish as $k\to \infty$. 

In the last section, we complete the proof of our main theorem. We find $\varepsilon>0$ depending on $\varepsilon_0$ and define a finite set of points which contains $\mathcal{S}$ and $\mathcal{S}'$:
$$
\mathcal{S}_0=\{x\in M: \lim_{r\to 0}\liminf_{k\to \infty}\int_{B_r(x)} R_k^2dV_{g_k}>\varepsilon\}.
$$ 
The points in $\mathcal{S}_0$ are said to be the bubble points. We apply the method of induction to show that if there is no real bubble at the bubble points, then the volume and the diameter of the bubble region will go to $0$. Moreover, since $\lambda_1(\Delta_{g_k})\geq\Lambda>0$, the sequence of manifolds $\{(M,g_k)\}$ can not be pinched as $k\to \infty$. Finally, the proof of the main theorem can be divided into 2 cases.  

\textbf{Case 1}: If $u\neq 0$, there is no real bubble at any bubble point, hence $\{(M,d_k)\}$ converges to $(M,d_u)$ in the Gromov-Hausdorff topology according to Proposition \ref{equal}. (See Figure1). 
\begin{figure}[htbp]
\centering
\begin{tikzpicture}[line cap=round,line join=round,>=triangle 45,x=0.7184358925374753cm,y=0.7669501590406798cm]
%\clip(0.6035726826804702,0.7461628741089976) rectangle (14.809928869370165,9.328018119181264);
\begin{scope}[scale=.75,xshift=0cm]
\draw [shift={(4.27518833306454,5.431695004496686)},line width=0.5pt]  plot[domain=1.5559897067437176:3.593755299774708,variable=\t]({1.*2.4690514198130393*cos(\t r)+0.*2.4690514198130393*sin(\t r)},{0.*2.4690514198130393*cos(\t r)+1.*2.4690514198130393*sin(\t r)});
\draw [shift={(0.8608343527345694,3.9718657416972745)},line width=0.5pt]  plot[domain=-0.32994878516225956:0.30907445871760614,variable=\t]({1.*1.2527945975832748*cos(\t r)+0.*1.2527945975832748*sin(\t r)},{0.*1.2527945975832748*cos(\t r)+1.*1.2527945975832748*sin(\t r)});
\draw [shift={(2.3770158234513703,3.254179925852358)},line width=0.5pt]  plot[domain=2.386021115432336:5.488988298291824,variable=\t]({1.*0.4546959047968492*cos(\t r)+0.*0.4546959047968492*sin(\t r)},{0.*0.4546959047968492*cos(\t r)+1.*0.4546959047968492*sin(\t r)});
\draw [shift={(3.5205869187532275,1.4120459050443472)},line width=0.5pt]  plot[domain=1.1839693752785962:2.0686202888026304,variable=\t]({1.*1.727473204911288*cos(\t r)+0.*1.727473204911288*sin(\t r)},{0.*1.727473204911288*cos(\t r)+1.*1.727473204911288*sin(\t r)});
\draw [shift={(4.282202803168843,5.472854399752712)},line width=0.5pt]  plot[domain=-1.615433323982165:0.5462446771001729,variable=\t]({1.*2.4634307134609097*cos(\t r)+0.*2.4634307134609097*sin(\t r)},{0.*2.4634307134609097*cos(\t r)+1.*2.4634307134609097*sin(\t r)});
\draw [shift={(6.992274478227374,7.3785864015371)},line width=0.5pt]  plot[domain=2.544373242953365:3.9439717690628098,variable=\t]({1.*0.8706729770525332*cos(\t r)+0.*0.8706729770525332*sin(\t r)},{0.*0.8706729770525332*cos(\t r)+1.*0.8706729770525332*sin(\t r)});
\draw [shift={(5.7173959052221806,7.873857327077636)},line width=0.5pt]  plot[domain=-0.010184963980616857:2.3670973291393795,variable=\t]({1.*0.5549457284211758*cos(\t r)+0.*0.5549457284211758*sin(\t r)},{0.*0.5549457284211758*cos(\t r)+1.*0.5549457284211758*sin(\t r)});
\draw [shift={(4.427723565648931,9.165657550849692)},line width=0.5pt]  plot[domain=4.620975218455563:5.491840563427147,variable=\t]({1.*1.2704864732558085*cos(\t r)+0.*1.2704864732558085*sin(\t r)},{0.*1.2704864732558085*cos(\t r)+1.*1.2704864732558085*sin(\t r)});
\draw [rotate around={45.8814039965821:(3.433986072291187,4.906829695395285)},line width=0.5pt] (3.433986072291187,4.906829695395285) ellipse (0.3551837133963182cm and 0.24521453417538153cm);
\draw [rotate around={45.8814039965821:(3.4257828157062895,6.670529861146692)},line width=0.5pt] (3.4257828157062895,6.670529861146692) ellipse (0.3551837133963024cm and 0.24521453417537065cm);
\draw [rotate around={45.8814039965821:(5.443783935589311,5.554886965601617)},line width=0.5pt] (5.443783935589311,5.554886965601617) ellipse (0.3551837133963277cm and 0.24521453417538805cm);
\draw [line width=0.5pt] (11.341925420125442,5.6041065051109555) ellipse (1.8207047799527298cm and 1.9436526418785973cm);
\draw [rotate around={45.8814039965821:(10.538006274806193,4.890423182225499)},line width=0.5pt] (10.538006274806193,4.890423182225499) ellipse (0.3551837133964161cm and 0.2452145341754491cm);
\draw [rotate around={45.8814039965821:(10.529803018221315,6.654123347976931)},line width=0.5pt] (10.529803018221315,6.654123347976931) ellipse (0.35518371339626453cm and 0.2452145341753445cm);
\draw [rotate around={45.8814039965821:(12.547804138104317,5.538480452431833)},line width=0.5pt] (12.547804138104317,5.538480452431833) ellipse (0.3551837133963403cm and 0.24521453417539682cm);
\draw [->,line width=0.5pt] (7.289516667189628,5.472854399752712) -- (8.372346536395149,5.489260912922492);
\end{scope}
\end{tikzpicture}
\caption{\label{fig: }When $u\neq 0$, $\{(M,d_k)\}$ converges to $(M,d_u)$ in the Gromov-Hausdorff topology.}
\end{figure}
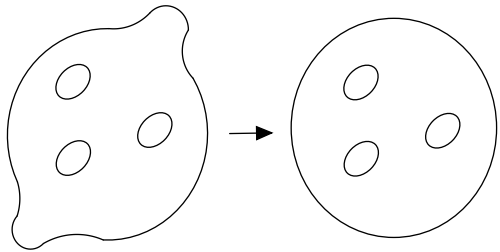

\textbf{Case 2}: If $u=0$, there is exactly one real bubble of $\{u_k\}$. After computation, we show $\liminf_{k\to \infty}\frac{1}{6}\int_{M} R_kdV_{g_k}\geq Y(\mathbb{S}^4)$ and get rid of this case. (See Figure2).
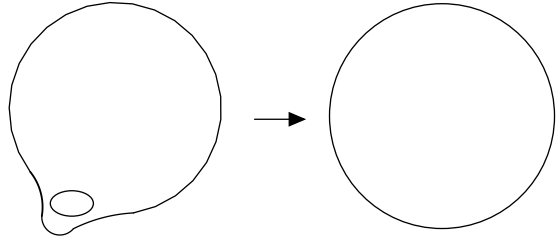
\begin{figure}[htbp]
\centering
\begin{tikzpicture}[line cap=round,line join=round,>=triangle 45,x=0.7cm,y=0.7cm]
%\clip(1.0728363895551223,0.6780106512175748) rectangle (12.88597129841931,6.537542898875724);
\draw [shift={(3.794848613394136,3.2852606182724426)},line width=0.5pt]  plot[domain=-1.398659944543927:3.782760967031551,variable=\t]({1.*2.0150407445358107*cos(\t r)+0.*2.0150407445358107*sin(\t r)},{0.*2.0150407445358107*cos(\t r)+1.*2.0150407445358107*sin(\t r)});
\draw [shift={(1.345,1.4125)},line width=0.7pt]  plot[domain=-0.16207294701362507:0.6743757121448317,variable=\t]({1.*1.069009471426705*cos(\t r)+0.*1.069009471426705*sin(\t r)},{0.*1.069009471426705*cos(\t r)+1.*1.069009471426705*sin(\t r)});
\draw [shift={(2.740322580645161,1.2208064516129034)},line width=0.5pt]  plot[domain=3.0852542314833955:5.5785109750612545,variable=\t]({1.*0.34086339081319683*cos(\t r)+0.*0.34086339081319683*sin(\t r)},{0.*0.34086339081319683*cos(\t r)+1.*0.34086339081319683*sin(\t r)});
\draw [shift={(4.293728813559325,-1.6001694915254299)},line width=0.5pt]  plot[domain=1.623753603150837:2.0324864803811344,variable=\t]({1.*2.904240972731641*cos(\t r)+0.*2.904240972731641*sin(\t r)},{0.*2.904240972731641*cos(\t r)+1.*2.904240972731641*sin(\t r)});
\draw [line width=0.5pt] (10.,3.14) circle (1.4956430055330718cm);
\draw [rotate around={0.:(2.97,1.48)},line width=0.5pt] (2.97,1.48) ellipse (0.2863790577818088cm and 0.1692689124913286cm);
\draw [->,line width=0.5pt] (6.44,3.08) -- (7.42,3.08);
\end{tikzpicture}
\caption{\label{fig: }When $u=0$, there is exactly one real bubble of $\{u_k\}$.}
\end{figure}
Thus, we finish proving our main theorem.

\textbf{Acknowledgements} The author is grateful to her doctoral supervisor Professor Yuxiang Li for his invaluable discussions and comments. 
\section{Preliminary}
Let $(M,g_0)$ be a compact 4-dim Riemannian manifold without boundary. Assume $\{g_k=u_k^2g_0\}$ is a sequence of conformal metrics and $u_k$ satisfies the following equation:
$$
-6\Delta u_k+R(g_0)u_k=R_ku_k^{3}.
$$
\subsection{$\varepsilon$-regulairity}
The goal of this subsection is to show the $\varepsilon$-regularity of conformal factors. Actually, in the foregoing paper \cite{Dong-Li-Xu}, we get the same regularity result of $\{u_k\}$ by using Moser Iteration. In this paper, we use another way to show the regularity of $\{u_k\}$ by modifying the definition of John-Nirenberg radius in \cite{LGZ}.
We consider the following operator over $\Omega\subset\mathbb{R}^4$:
$$
Lu=-div(a^{ij}u_j)+cu,
$$
where
\begin{align}\label{2}
0<\lambda_1\leq a^{ij}, \s\s ||\nabla a^{ij}||_{C^0(\Omega)}+||a^{ij}||_{C^0(\Omega)}<\lambda_2, \s\s ||c||_{C^0(\Omega)}<\lambda_3.
\end{align}
 Let $\{u_k\}$ be a sequence of positive functions in $W^{1,2}(\Omega)$, each of which solves the equation $Lu_k=f_ku_k$ in the weak sense, where $||f_k||_{L^2(\Omega)}$ is uniformly bounded. 

Now, recall the following John Nirenberg lemma, which was presented in \cite{Moser}.
\begin{lem}
If $w$ is a square integrable in the unit cube $Q(1)\subset \R^n$ and if for every cube $Q\subset Q(1)$
\begin{equation*}
\frac{1}{|Q|}\int_Q(w-w_Q)^2dx\leq 1, where, w_Q=\frac{1}{|Q|}\int_Qwdx,
\end{equation*}
then every power of $w$ is integrable and even more: There exist positive constants $\alpha, \beta$ depending on $n$ only, such that
\begin{equation*}
\int_{Q(1)}e^{\alpha w}\int_{Q(1)}e^{-\alpha w}\leq \beta^2.
\end{equation*}
\end{lem}
Replacing the cube $Q(1)$ with the inscribed ball $B_r(x)$, we define John-Nirenberg radius as follows.
\begin{defi}\label{radius}
Given $u$ is a positive function, we define John-Nirenberg radius, for any $x\in \Omega_1\subset\subset\Omega\subset \mathbb{R}^4$ and $\varepsilon>0$ 
\begin{center}
$\rho(x,u,\Omega,\varepsilon)=\sup\{r:\frac{1}{|B_{t}(x)|}\int_{B_{t}(x)}|\log u-\overline{(\log u)}_{B_t(x)}|<\varepsilon, \forall t\leq r,  B_{r}(x)\subset \subset \Omega\}$,
\end{center}
where  $\overline{(\log u)}_{B_t(x)}$ is the mean integral of  $\log u$ over $B_t(x)$ .
\end{defi}
This kind of radius is a crucial point to study the convergence of $\{u_k\}$ in this subsection. The following proposition is given at first. 
\begin{pro}\label{positiveradius}
For any $\varepsilon>0$, there exist $\varepsilon_0>0$ and $a>0$ which only depend on $\lambda_1,\lambda_2,\lambda_3$ and $\varepsilon$, such that if 
$$\int_{B_3} |f_k|^{2}dx<\varepsilon_0, $$
then 
$$
\rho(x,u_k,B_3,\varepsilon)>a>0,  \forall x\in B_1\subset B_3,
$$
when $k$ is sufficiently large. $B_r$ denotes a ball of radius $r$ centered at $0$ in $\R^4$.
\end{pro}
The following lemma  can be derived from the above proposition and John-Nirenberg lemma immediately.
\begin{lem}\label{regularity0}
Let $\Omega$ be a domain of $\mathbb{R}^4$ and $\int_{\Omega} |f_k|^2dV_{g_k} $ be uniformly bounded for any $k$. For any $p\in(1,2)$, there exits $\varepsilon>0$, such that if for any $x\in \Omega_1\subset\subset \Omega$, $\rho(x,u_k,\Omega,\varepsilon)> a > 0$, then we can find $\{c_k\}$ such that $\{c_ku_k\}$, $\{\log (c_ku_k)\}$ and $\{\frac{1}{c_ku_k}\}$ converge to $u$, $\log u$ and $\frac{1}{u}$ weakly in $W^{2,p}(\Omega_2)$, where $\Omega_2\subset\subset\Omega_1$.
\end{lem}
\proof
Take $\phi=\eta^2u_k^{-1}$ as a test function. For $\Omega_1\subset\subset\Omega$, take $\eta\equiv 1$ on $\Omega_1$, $|\nabla \eta|<C$ and $\eta\in C^\infty_0(\Omega)$. Multiplying $Lu_k=f_ku_k$ by $\phi$ and integrating, we have
$$
\int_\Omega \left(a^{ij}\nabla_j u_k\nabla_i(\frac{\eta^2}{u_k})+c\eta^2\right)dx=\int_\Omega f_k\eta^2.
$$
Using the Young's inequality and the H\"older's inequality obtains
$$
\int_\Omega\eta^2|\nabla \log u_k|^2\leq C_1+C_2(\int_\Omega f_k^2)^{\frac{1}{2}},
$$
where $C_1$ and $C_2$ depend on $\lambda_1,\lambda_2$ and $\lambda_3$. Hence $\nabla \log u_k$ has an uniform $L^2$-norm upper bound over $\Omega_1$. Choose $c_k$ such that $\int_{\Omega_1}\log c_ku_k=0$, by the Poincar\'e inequality, $||\log c_ku_k||_{W^{1,2}(\Omega_1)}$ are uniformly bounded. According to the Sobolev embedding theorem, $\{\log c_ku_k\}$ converges to a function $\varphi$ weakly in $L^{4}$.

Now, take finite open balls $B_a(x_1),\cdots,B_a(x_n)$ cover $\bar{\Omega}_1$. In each ball $B_a(x_i)$, there exists $c_k^i$ such that $\{\log c_k^iu_k\}$ converges to  $\varphi^i$  pointwise a.e. $x\in B_a(x_i)$. Using the Egorov's theorem, we derive $\{\log c_k^iu_k\}$ uniformly converges to $\varphi^i$ on a measurable set $E_i\subset B_a(x_i)$ with $\mathcal{L}(B_a(x_i)\setminus E_i)<\varepsilon'$, where $\varepsilon'$ is sufficiently small and $\mathcal{L}$ is the Lebegue measure over $\mathbb{R}^4$.  Since $\rho(x_i,u_k,\Omega,\varepsilon_0)> a > 0$, we get
$$
\int_{B_a(x_i)} e^{\frac{\alpha}{\varepsilon}\log c_k^iu_k}\int_{B_a(x_i)} e^{-\frac{\alpha}{\varepsilon}\log c_k^iu_k}\leq \beta,
$$
where $\alpha$ and $\beta$ both depend on the dimension. Since $\{e^{\frac{\alpha}{\varepsilon}\log c_k^iu_k}\}$ and $ \{e^{-\frac{\alpha}{\varepsilon}\log c_k^iu_k}\}$ converge to $e^{\frac{\alpha}{\varepsilon}\varphi^i}$ and $e^{-\frac{\alpha}{\varepsilon}\varphi^i}$ respectively over $E_i$, we have
$$
\int_{E_i} e^{\frac{\alpha}{\varepsilon}\log c_k^iu_k}\geq C_1>0 \s and \s
\int_{E_i} e^{-\frac{\alpha}{\varepsilon}\log c_k^iu_k}\geq C_2>0,
$$
which imply 
$$
\int_{B_a(x_i)} e^{\frac{\alpha}{\varepsilon}\log c_k^iu_k}<C \s and \s \int_{B_a(x_i)} e^{-\frac{\alpha}{\varepsilon}\log c_k^iu_k}<C.
$$
We may assume $\{\log c_k-\log c_k^i\}$ converges, then
$$
\int_{B_a(x_i)} e^{\frac{\alpha}{\varepsilon}\log c_ku_k}<C \s and \s \int_{B_a(x_i)} e^{-\frac{\alpha}{\varepsilon}\log c_ku_k}<C.
$$
It turns out that both $||c_ku_k||_{L^{\frac{\alpha}{\varepsilon}}(\Omega_1)}$ and $||\frac{1}{c_ku_k}||_{L^{\frac{\alpha}{\varepsilon}}(\Omega_1)}$ are bounded.

Take $w_k=c_ku_k$ and consider the following equation:
$$
-div(a^{ij}(w_k)_j)+cw_k=f_kw_k.
$$
For any $p\in(1,2)$, take $\varepsilon$ small enough, such that $p\in (1, \frac{2\alpha}{2\varepsilon+\alpha}]\subset (1,2)$,
$$
\int_{\Omega_1}(f_kw_k)^p\leq (\int_{\Omega_1}f_k^2)^\frac{p}{2} (\int_{\Omega_1} w_k^{\frac{2p}{2-p}})^\frac{2-p}{2}<C.
$$
Then we get $||w_k||_{W^{2,p}(\Omega_2)}$ is bounded for any $\Omega_2\subset\subset \Omega_1$, by the standard elliptic theory. Similarly, the estimate of $||w_k^{-1}||_{W^{2,p}(\Omega_2)}$ and $||\log w_k||_{W^{2,p}(\Omega_2)}$ can be derived from the following equalities:
$$
\nabla w_k^{-1}=-\frac{\nabla w_k}{w_k^2}, \s \nabla^2 w_k^{-1}=-\frac{\nabla^2 w_k}{w_k^2}+2\frac{|\nabla w_k|^2}{w_k^3},
$$
and
$$
\nabla\log w_k=\frac{\nabla w_k}{w_k}, \s \nabla^2 \log w_k = \frac{\nabla^2 w_k}{w_k}-2\frac{|\nabla w_k|^2}{w_k^2}.
$$
\endproof
\begin{rem}
In the above proof, the underlying meaning of $\{c_k\}$ and $\{c_k^i\}$ is to make the Poincar\'e inequality hold. Hence, there are many different ways to choose $\{c_k\}$ and $\{c_k^i\}$.
\end{rem}

\textbf{proof of Proposition \ref{positiveradius}}

Assume the result is not true. We can find a subsequence of functions $\{u_k\}$ and a sequence of points  $\{x_k\}\subset B_1$ such that
 $$
 \int_{B_3} |f_k|^2dx\to 0,
 $$
 and 
 \begin{equation*}
 \rho(x_k,u_k,B_3,\varepsilon)=\inf_{x\in B_1}\rho(x,u_k,B_3,\varepsilon) \to 0.
 \end{equation*}
Denote $\rho_k(x)$ as $\rho(x,u_k,B_3,\varepsilon)$
and find $y_k \in B_2$ such that
\begin{align*}
\frac{\rho_k(y_k)}{2-|y_k|}&=\inf_{x\in B_2}\frac{\rho_k(x)}{2-|x|}\\&\leq \frac{\rho_k(x_k)}{2-|x_k|}\\&\leq\rho_k(x_k)\to 0. 
\end{align*}
Thus, $\rho_k(y_k)\to 0$. For any fixed $R$ and $y\in B_{R\rho_k(y_k)}(y_k)\subset B_{2-|y_k|}(y_k)\subset B_2$, it holds
\begin{align*}
\frac{\rho_k(y)}{\rho_k(y_k)} &\geq\frac{2-|y|}{2-|y_k|}\\ &\geq 1-\frac{R\rho_k(y_k)}{2-|y_k|}\\ &> \frac{1}{2},
\end{align*}
when $k$ is sufficiently large.

Now, take a subsequence $\{y_k\}$ which converges to $y_0$ and set $v_k(x)=r_ku_k(y_k+r_k x)$ where $r_k=\rho_k(y_k)$. Then we have
\begin{center}
$\rho(x,v_k,B_R,\varepsilon)> \frac{1}{2}$,  $\forall x\in B_{\frac{R}{2}}$,
\end{center}
and
$$
\frac{1}{|B_1|}\int_{B_1}|\log v_k-\overline{(\log v_k)}_{B_1(x)}|=\varepsilon.
$$
 $v_k$ satisfies $-div(a(x_k+r_kx)^{ij}(v_k)_j)+c(x_k+r_kx)r_k^2v_k=r_k^2f_k(y_k+r_kx)v_k$. Notice
  $$
\int |f_k(x)|^2r_k^4(x)dx=\int |f_k(y_k+r_k x)|^2d(x_k+r_k x )\to 0.
$$
By the above lemma, there exists $\{c_k\}$,  such that $\frac{1}{|B_1|}\int_{B_1}|\log c_kv_k-\overline{(\log c_kv_k)}_{B_1(x)}|=\varepsilon$ and $\int_{B_1} \log c_ku_k=0$, then $\{c_kv_k\}$ weakly converges to a positive harmonic function $v_0$ satisfying $\int \log v_0=0$ in $W^{2,p}(\R^4)$.  By the Liouville's theorem, $v_0$ is a constant. It is in contradiction to $\frac{1}{|B_1|}\int_{B_1}|\log v_0-\overline{(\log v_0)}_{B_1(x)}|=\varepsilon$. Thus we finish the proof.
\endproof

Now, the $\varepsilon$-regularity of $\{u_k\}$ is given as follows.
\begin{lem}\label{regularity1}
\textbf{($\varepsilon$-Regularity)} For any $p\in(1,2)$, there exists $\varepsilon_0 >0$ depending on $p$, such that if $$\int_{B_r}|f_k|^2 dx<\varepsilon_0,$$ then
$$
||u_k||_{W^{2,p}(B_{\frac{r}{2}})}\leq C ||u_k||_{L^4(B_r)},
$$
where $C$ only depends on $\lambda_1$, $\lambda_2$ and $\lambda_3$ in (\ref{2}).
\end{lem}
\proof
Assume the result fails to hold. We can find a subsequence of functions $\{u_k\}$, such that 
$$\int_{B_r(x_0)}|f_k|^2 dx\to 0,$$ 
and
$$
||u_k||_{W^{2,p}(B_{\frac{r}{2}})}\geq k||u_k||_{L^4(B_r)}.
$$
From Proposition \ref{positiveradius}, if $\int_{B_r}|f_k|^2 d\mu_{g_k}\to 0$, then $\rho(x,u_k,B_r,\varepsilon)>a>0, \s\forall x\in B_{\frac{3}{4}r}$. Then take $\varepsilon$ be sufficiently small, by Lemma \ref{regularity0}, we find $\{c_k\}$ such that $\{c_ku_k\}$ weakly converges to $ v$ in $W^{2,p}(B_{\frac{3}{5}r})$.
Besides, $||c_ku_k||_{L^\frac{\alpha}{\varepsilon}(B_{\frac{3}{5}r})}$ and $||\frac{1}{c_ku_k}||_{L^\frac{\alpha}{\varepsilon}(B_{\frac{3}{5}r})}$ are bounded. By the H\"older's inequality
\begin{equation*}
0<\int_{B_{\frac{3}{5}r}} dx\leq \left(\int_{B_{\frac{3}{5}r}}|c_ku_k|^4dx\right)^\frac{1}{2}\left(\int_{B_{\frac{3}{5}r}}|\frac{1}{c_ku_k}|^4dx\right)^\frac{1}{2},
\end{equation*}
hence $||v||_{L^4(B_{\frac{3}{5}r})}$ is positive. Thus
$$
||c_ku_k||_{W^{2,p}(B_{\frac{r}{2}})}\geq  k||c_ku_k||_{L^4({B_{\frac{3}{5}r}}) }\to \infty,
$$
which is impossible.
\endproof

Now, assume $\{g_k=u_k^2g_0\}$ is a sequence of conformal metrics over $(M,g_0)$. And apply Lemma \ref{regularity1} to a small ball $B_r^{g_0}(x)\subset (M,g_0)$ with
$$
f_k=R_ku_k^2 \s and\s c=R_0.
$$
Then we get the following $\varepsilon$-regularity of $u_k$ in $(M,g_0)$.
\begin{lem}\label{regularity}
Let $B_r(x_0)\subset (M,g_0)$. For a sequence of conformal metrics $g_k=u_k^2g_0$ in the fixed conformal class $[g_0]$ and $p\in(1,2)$, there exists $\varepsilon_0>0$, such that if 
$$
\int_{B_r(x_0)} R_k^2dV_{g_k}<\varepsilon_0,
$$
then 
$$
||u_k||_{W^{2,p}(B_\frac{r}{2}(x_0))}<C \mathrm{Vol}(B_r(x_0),g_k)^\frac{1}{4},
$$
where $C$ only depends on $(M,g_0)$.
\end{lem}

\subsection{Removability of singularities} In this subsection, we  focus on the set of concentration points about $\{u_k\}$:
$$
\mathcal{S}=\{x\in M:\lim_{r\to 0}\liminf_{k\to \infty} \int_{B_r(x)} R_k^2dV_{g_k}> \varepsilon_0\}.
$$
It is easy to show $\mathcal{S}$ is a finite set. According  to Lemma \ref{regularity}, for $p<2$, $\{u_k\}$ converges to $u$ weakly in $W^{2,p}$ on any compact set $K\subset M\setminus \mathcal{S}$. If $\mathcal{S}$ is nonempty, it's possible that $u$ has a singularity at some $x\in \mathcal{S}$, but in fact this is not the case. 

At first, we establish an important tool Three Circles Lemma. Let $Q=[0,3L]\times \mathbb{S}^{3}$ and $Q_i=[(i-1)L,iL]\times \mathbb{S}^{3}, i=1,2,3$. Set $g_{q}=dt^2+g_{\mathbb{S}^{3}}$ and $dV_q=dV_{g_{q}}$.
\begin{lem}\label{TCL}
Assume $g$ is a Riemannian metric over $Q$ and $u\in W^{1,2}(Q)$ is the weak solution of equation $-6\Delta u+R_gu=Ru^3$, where $R_g$ and  $R$ are the scalar curvature of $g$ and $u^2g$ respectively. If there exist $\varepsilon_1>0$, $\tau_1>0$ and $L_0>0$ such that 
$$||g-g_q||_{C^2}<\tau_1 \s and \s \int_Q R^2u^4dV_g<\varepsilon_1,$$
then for any $L>L_0$ we have

(1)$\int_{Q_1} u^2dV_q\leq e^{-L}\int_{Q_2} u^2dV_q$ implies
$$
\int_{Q_2} u^2dV_q\leq e^{-L}\int_{Q_3} u^2dV_q.
$$

(2)$\int_{Q_2} u^2dV_q\geq e^{-L}\int_{Q_3} u^2dV_q$ implies
$$
\int_{Q_1} u^2dV_q\geq e^{-L}\int_{Q_2} u^2dV_q.
$$

(3)either $\int_{Q_1} u^2dV_q\leq e^{-L}\int_{Q_2} u^2dV_q$ or $\int_{Q_2} u^2dV_q\leq e^{-L}\int_{Q_1} u^2dV_q$.
\end{lem}
 We omit the proof, since it is very similar to \cite{Li-Zhou}(cf. Theorem 4.1). 

Set $L_{g_0}=-6\Delta+R_0$. Assume $\{R_ku_k^2\}$ converges to $f$ in distribution sense over $M$. We check the following lemma.
\begin{lem}\label{definescalar}
$u$ is a weak solution of the following elliptic equation over $M\setminus \mathcal{S}$:
$$L_{g_0}u=fu.$$
\end{lem}
\proof
For simplicity, assume $\mathcal{S}=\{P\}$. Let $\phi\in C^\infty_0(M\setminus\{P\})$. $\{u_k\}$ denotes a subsequence which $W^{2,p}$ weakly converges to $u$ on the compact support of $\phi$. 
\begin{align*}
\int (L_{g_0}\phi)udV_0&=\int(L_{g_0}\phi)(u-u_k)dV_0+\int(L_{g_0}\phi)u_kdV_0
\\ &=\int(L_{g_0}\phi)(u-u_k)dV_0+\int\phi (L_{g_0}u_k)dV_0
\\&=\int(L_{g_0}\phi)(u-u_k)dV_0+\int\phi R_ku_k^2u_kdV_0
\\&=\int(L_{g_0}\phi)(u-u_k)dV_0+\int\phi R_ku_k^2(u_k-u)dV_0+\int\phi R_ku_k^2udV_0.
\end{align*}
Since $u_k$ $L^q$-converges to $u$ over the support of $\phi$ for some $1<q<\frac{4p}{4-2p}$. The first two terms go to $0$ as $k\to \infty$ and $f$ is a weak $L^{2}$-limit of $\{R_ku_k^2\}$, then  
$$
\int (L_{g_0}\phi)udV_0=\int\phi fudV_0.
$$
\endproof

Hence, we can define $R(g')=u^{-2}f$ as the scalar curvature of $g'=u^2g_0$ over $M$. By Lemma \ref{definescalar}, $u$ satisfies $-6\Delta u+R_0u=R(g')u^3$ in the weak sense on $M\setminus \mathcal{S}$, and $\int_M R(g')^2u^4dV_{g_0}=\int_M f^2dV_{g_0}$ is bounded. Then the removability of singularities can be derived from Lemma \ref{TCL}.
\begin{lem}
Let $g_0$ be a smooth metric over $B_1$ and $u\in W^{2,p}(B_1\setminus \{0\},g_0)$, for some $p\in(1,2)$. If 
$$
\int_{B_1}(1+R(g')^2)dV_{g'}<+\infty,
$$
then $g'$ can be extended in $W^{2,p}(B_1)$.
\end{lem}
\proof
Choose a normal chart around $0$ with respect to $g_0$. Set 
\begin{align*}
\phi(r,\theta)=(-\log r,\theta),
\end{align*} 
and
\begin{align*}
g'(t,\theta)=v^2\hat{g}(t,\theta),
\end{align*}
where $\hat{g}(t,\theta)=\phi^*(g_0)(t,\theta)$. Then
\begin{align*}
v(t,\theta)=u(e^{-t},\theta)e^{-t}, \s and \s  \int_{B_1}(1+R(g')^2)u^4dx=\int_{[0,+\infty]\times\mathbb{S}^3} (1+R(g')^2)v^4dV_q.
\end{align*}
Given a sufficiently small $\delta>0$, we assume
\begin{align*}
||\hat{g}-g_{q}||_{C^2([-\log \delta,\infty]\times \mathbb{S}^{3})}\leq \tau_1,\s \s  \int_{[-\log \delta,\infty]\times \mathbb{S}^{3}}|R(g'(v))|^2dV_{g'(v)}<\min\{\varepsilon_0, \varepsilon_1\}.
\end{align*}
where $\varepsilon_0$ $\varepsilon_1$ and $\tau_1$ are mentioned in Lemma \ref{regularity} and Lemma \ref{TCL}. $C$ denotes different bounded constants. Take $L\geq L_0$, according to Lemma \ref{TCL}, for any $i\geq 0$, we have
\begin{align*}
\int_{[-\log \delta+iL,-\log \delta+(i+1)L]\times \mathbb{S}^{3}} v^2dV_q  &\leq e^{-iL}\max\{\int_{[-\log\delta,-\log\delta+L]\times \mathbb{S}^{3}}v^2dV_q, 
\\ &\int_{[-\log\delta+2iL,-\log\delta+(2i+1)L]\times \mathbb{S}^{3}}v^2dV_q\}
\\&\leq e^{-iL} ||v||_{L^2([-\log \delta,\infty]\times \mathbb{S}^{3})}^2.
\end{align*}
Then 
\begin{align*}
\int_{[-\log \delta+iL,-\log \delta+(i+3)L]\times \mathbb{S}^{3}} v^{4}dV_q&\leq C(\int_{[-\log \delta+iL,-\log \delta+(i+3)L]\times \mathbb{S}^{3}} v^2dV_q)^\frac{1}{2}
\\ &\leq C(e^{-iL}+e^{-(i+1)L}+e^{-(i+2)L})^\frac{1}{2}
\\ & \leq Ce^{i\frac{L}{2}}.
\end{align*}
Hence
\begin{align*}
\left(\int_{B_{e^{-iL}\delta}\setminus B_{e^{-(i+3)L}\delta}} u^{4}dx\right)^\frac{p}{4}= \left(\int_{(-\log \delta+iL,-\log \delta+(i+3)L)\times \mathbb{S}^{3}} v^4 dV_q\right)^\frac{p}{4}<Ce^{-i\frac{Lp}{8}}.
\end{align*}
By $\varepsilon$-regularity: 
$$
||u||_{W^{2,p}(B_{e^{-iL}\delta}\setminus B_{e^{-(i+1)L}\delta})}\leq C||u||_{L^4(B_{e^{(-i+1)L}\delta}\setminus B_{e^{-(i+2)L}\delta})}.
$$
Hence
$$
||u||^p_{W^{2,p}(B_{e^{-L}\delta})}=\sum_{i=1}^\infty||u||^p_{W^{2,p}(B_{e^{-iL}\delta}\setminus B_{e^{-(i+1)L}\delta})}
\leq C\sum_{i=0}^\infty||u||^p_{L^4(B_{e^{-iL}\delta}\setminus B_{e^{-(i+3)L}\delta})}<C\sum_{i=0}^\infty e^{-i\frac{Lp}{8}}<\infty.
$$
\endproof

\section{Local Convergence of Distance Functions}
It's known that a connected Riemannian manifold carries the structure of a metric space with the distance function defined by the arc length of a minimizing geodesic. In this section, we
focus on the distance functions $d_k$ defined by $g_k=u_k^2g_0$ over $M\times M$.
\begin{align*}
d_k(x,y)&=\inf_{\mbox{piecewise smooth $\gamma\subset M$ joining $x$ and $y$}}\int_\gamma \sqrt{g_k(\dot{\gamma},\dot{\gamma})}
\\ &= \inf_{\mbox{piecewise smooth $\gamma\subset M$ joining $x$ and $y$}}\int_\gamma u_k.
\end{align*}

In order to show local convergence of $\{d_k\}$, we present the following series of lemmas. Denote the canonical Euclidean metric over Euclidean space $\R^4$ by $g_{eucl}$.
Let $g=u^2g_{eucl}$ be a conformal metric defined on $\R^4$ and $\gamma$ be a piecewise smooth curve connecting $\gamma(0)=0$, $\gamma(1)=y$. We define the distance function $d_u$ as follows:
\begin{align*}
d_u(0,y)= \inf_{\gamma\subset \R^4}\int_\gamma u.
\end{align*}

\begin{lem} \label{lem3}
For any $\varepsilon'>0$, there exist $\varepsilon_2>0,\s\tau_2>0$, which depend on $\varepsilon'$,  such that if 
$$\int_{B_3}R_u^2dV_u<\varepsilon_2,$$
then
$$
 (1+\varepsilon')d_u(0,y)\geq\inf_{\gamma\subset B_2}\int_\gamma u,\s \forall y\in B_{\tau_2},
$$
where $R_u$ denotes the scalar curvature of $u^2g_{eucl}$.
\end{lem}
\proof
If the result is not true. Given $\varepsilon'>0$, we can find a sequence of points $\{y_k\}\subset B_2$ and a sequence of conformal metrics $\{g_k=u_k^2g_{eucl}\}$ such that
$$
r_k=|y_k|\to 0,\s \int_{B_3}R_k^2dV_{g_k}\to 0, \s and \s (1+\varepsilon')d_k(0,y_k)<\inf_{\gamma\subset B_2}\int_\gamma u_k.
$$
Take $v_k=c_kr_ku_k(r_kx)$, where $c_k$ is chosen such that 
$$
0=\int_{B_1} \log v_k.
$$
Since for any fixed $R$:
$$
\int_{B_R} \hat{R}_k^2dV_{\hat{g}_k}\leq \int_{B_3}R_k^2dV_{g_k}\to 0, 
$$
where $\hat{R}_k$ is the scalar curvature with respect to $\hat{g}_k=v_k^2g_{eucl}$. By Lemma \ref{regularity}, we may assume $\{v_k\}$ weakly converges to a positive function $v$ in $W^{2,p}(\R^4)$ with $\int_{B_1} \log v_k=0$, for some $p<2$. By the Liouville's theorem, $v=1$. According to the Trace Embedding Theorem \cite{Adam} (cf. Theorem 4.12), we have  
$$
c_kr_k\inf_{\gamma\subset B_2}\int_\gamma u_k\leq \int_0^{y_k} c_kr_ku_k=\int_0^1 c_kr_ku_k(r_kx)\to 1.
$$
Using $\hat{d}_k$ denote the distance function defined by $\hat{g}_k$,  we get (see \cite{Dong-Li-Xu} Proposition 3.2):
$$
c_kr_kd_k(0,y_k)=\hat{d_k}(0,1)\to 1.
$$ 
Hence we get a contradiction. The proof is finished.
\endproof
Letting $\varepsilon'\to 0$, we get 
$$
d_u(0,y)\geq \inf_{\gamma\subset B_2}\int_{\gamma} u, \s \forall y\in B_{\tau_2}.
$$
According to the definition of $d_u$, it is obvious that 
$$
d_u(0,y)=\inf_{\gamma\subset B_2}\int_{\gamma} u, \s \forall y\in B_{\tau_2}. 
$$

Combining the above lemma and the main theorem which we proved in the previous paper \cite{Dong-Li-Xu} (cf. Theorem 1.2), we define a finite set of points:
$$
\mathcal{S}'=\{x\in M:\lim_{r\to 0}\liminf_{k\to 0}\int_{B_r(x)}R_k^2dV_{g_k}>\min\{\varepsilon_0,\varepsilon_2\}\}.
$$
With a simple covering argument, we can get the following lemma. 
\begin{lem} \label{lem4}
For any $x\notin \mathcal{S}'$, there exists $r>0$, such that $\{u_k\}$ converges to $u$ weakly in $W^{2,p}(B_r(x))$ and $d_k$ uniformly converges to $d_u$ on $B_r(x)$ for some $p<2$. If $u\neq 0$, $d_u$ is the distance function defined by $u^2g_0$ over entire $M\times M$.
\end{lem}
Since $u\in W^{2,p}(M,g_0)$, $d_u$ defined by $u^2g_0$ is continuous over entire $M\times M$. If $u\equiv 0$, we can find a sequence $\{c_k\}$, the above lemma remains valid for $\{v_k=c_ku_k\}$. We now derive the convergence of $\{(M,d_k)\}$ in the Gromov-Hausdorff topology as follows.
\begin{pro}\label{equal}
If
$$
\lim_{r\rightarrow 0}\mathrm{diam}(B_r(P),g_k)+\mathrm{Vol}(B_r(P),g_k)=0,\s \forall P\in\mathcal{S}',
$$
then  $(M,d_{k})$ converges to $(M,d_u)$ in the sense of 
Gromov-Hausdorff distance. 
\end{pro}
\proof
Without loss of generality, we assume $\mathcal{S}'=\{P\}$ and $R_k^2dV_{g_k}$
converges to $\mu$ in the sense of distribution. Set $\varepsilon=\min\{\varepsilon_0,\varepsilon_2\}$.

Given a $r>0$, we can choose $r'<\frac{r}{4}$, such that 
$\mu(B_{2r'}(x))<\frac{\varepsilon}{2}$, for any $x\notin B_r(P)$.
Then $\int_{B_{r'}(x)}R_k^2dV_{g_k}<\varepsilon$ for any $x\in M\setminus B_r(P)$, when $k$ is sufficiently
large. Then we may assume $d_{k}$ induced by $g_k$ converges to a
function $d_\infty$ on $C^0(M\setminus B_r(p))$.

For any $x,y\in B_r(P)\setminus\{P\}$,  
$$
d_\infty(x,y)=\lim_{k\rightarrow+\infty}d_{k}(x,y)\leq
\lim_{k\rightarrow+\infty}\mathrm{diam}(B_r(P),g_k).
$$
Thus $\lim_{x\rightarrow P}d_\infty(x,y)$ exists for any $y$.
We define
$$
d_\infty(P,y)=\lim_{x\rightarrow P}d_\infty(x,y).
$$
It is easy to check $d_\infty$ is a distance function and
$$
(M,d_{k})\xrightarrow{G-H}(M,d_\infty).
$$
Since $\lim_{r\rightarrow 0}\mathrm{Vol}(B_r(P),g_k)=0$ and $\mathrm{Vol}(M,g_k)=1$, we get $u\neq 0$. Hence $d_u$ is the distance function defined by $u^2g_0$.

Next, we prove $d_\infty=d_u$. Given $x,y\neq P$,

1)$d_\infty(x,y)\leq d_u(x,y)$.

For any piecewise smooth curve $\gamma$ connecting $x$ and $y$. Assume $\gamma(0)=x$ and $\gamma(1)=y$. Choose $t_1$ and $t_2$, such that $\gamma(t_i)\in \partial B_r(P)(i=1,2)$ and $\gamma|_{[0,t_1]\cup[t_2,1]}\cap B_r(P)=\emptyset$. Then 
\begin{align*}
d_\infty(x,y)&\leq d_\infty(x,\gamma(t_1))+d_\infty(\gamma(t_2),y)+\mathrm{diam}(B_r(P),d_\infty)
\\ &= \lim_{k\to \infty} d_k(x,\gamma(t_1))+\lim_{k\to \infty}d_k(\gamma(t_2),y)+\mathrm{diam}(B_r(P),d_\infty)
\\ &\leq \lim_{k\to \infty} \int_{\gamma|_{[0,t_1]}} u_k +\lim_{k\to\infty}\int_{\gamma|_{[t_2,1]}} u_k+\mathrm{diam}(B_r(P),d_\infty)
\\ &=\int_{\gamma|_{[0,t_1]}} u +\int_{\gamma|_{[t_2,1]}} u+\mathrm{diam}(B_r(P),d_\infty),
\end{align*}
the last equality comes from the Trace Embedding Theorem \cite{Adam}. Letting $r\to 0$, we get $\mathrm{diam}(B_r(P),d_\infty)=\lim_{k\to 0}\mathrm{diam}(B_r(P),d_k)=0$ and 
$$
d_\infty(x,y) \leq \int_{\gamma|_{[0,t_1]}} u +\int_{\gamma|_{[t_2,1]}} u\leq \int_\gamma u.
$$
Hence $d_\infty(x,y)\leq d_u(x,y)$.

2)$d_\infty(x,y)\geq d_u(x,y)$.

We set $\gamma$ be the segment defined in $(M,d_\infty)$ connecting $x$ and $y$, i.e. $\gamma: [0,L]\to (M,d_\infty)$ is a continuous map which satisifies
$$
d_\infty(\gamma(s),\gamma(s'))=|s-s'|, \s \forall s,s'\in [0,L].
$$

Consider the following two cases:

First, we consider the case when $d_0(P,\gamma)=\delta>0$. We claim that $\gamma$ is also continuous in $(M,g_0)$. Assume this is not true, there exists $t_i\to t$ and $\rho>0$ such that $d_0(\gamma(t_i),\gamma(t))>\rho$. Since whenever $d_0(\gamma(t_i),\gamma(t))>\rho$, there exists $\tau>0$, which depends on $g_0$ and $||u_k^{-1}||_{W^{1,q}} (q<4)$, such that $d_k(\gamma(t_i),\gamma(t))>\tau$ (\cite{Dong-Li-Xu} c.f.(3.1)). And $d_k$ uniformly converges to $d_\infty$, hence 
$$|t_i-t|=\lim_{k\to \infty} d_k(\gamma(t_i),\gamma(t))>\tau, $$
which is impossible. By Lemma \ref{lem4}, $d_u(x,y)=d_\infty(x,y)$
when $d_{0}(y,x)<\frac{\tau_2}{3}d_{0}(P,x)$. Hence we can find finitely many
points $p_0=\gamma(0)=x$, $p_1=\gamma(t_1)$, $\cdots$, $p_m=\gamma(1)=y$, such that $\forall i=0,\cdots, m-1$:
$$
d_0(p_i,p_{i+1})<\frac{\tau_2}{3} min\{d_0(p_0,P),\cdots, d_0(p_m,P)\}\leq \frac{\tau_2}{3} \delta.
$$
Thus
$d_\infty(p_i,p_{i+1})=d_u(p_i,p_{i+1})$ and 
$d_\infty(x,y)=\sum_id_u(p_i,p_{i+1})\geq d_u(x,y)$.

Then we focus on the case in which $d_0(P,\gamma)=0$. 
For any fixed $r>0$, choose $t_1$ and $t_2$, such that $\gamma(t_i)\in \partial B_r(P)(i=1,2)$ and $\gamma|_{[0,t_1]\cup[t_2,1]}\cap B_r(P)=\emptyset$. Then 
\begin{align*}
d_u(x,y)&\leq d_u(x,\gamma(t_1))+d_u(\gamma(t_2),y)+d_u(\gamma(t_1),\gamma(t_2))
\\ &= d_\infty(x,\gamma(t_1))+d_\infty(\gamma(t_2),y)+d_u(\gamma(t_1),\gamma(t_2))
\\ &\leq  d_\infty(x,y)+\mathrm{diam}(B_r(P), u^2g_0).
\end{align*}
Since $u\in W^{2,p}(B_r(P))$, we get 
$$
\lim_{r\to 0}\mathrm{diam}(B_r(P), u^2g_0)=0.
$$
Thus $d_\infty(x,y)\geq d_u(x,y)$.
\endproof

\section{Neck Analysis}
At first, we recall some definitions about bubble trees which will be used in the following analysis. 
$\{(x_k,r_k)\}$ is called a nontrival blowup sequence at $x$, if $r_k\rightarrow 0$, $x_k\rightarrow x$ and there exists a finite
set $\mathcal{S}$, such that 
$\{r_ku_k(x_k+r_kx)\}$ converges weakly in $W^{2,p}_{loc}(\R^4\setminus
\mathcal{S})$ to a positive function for some $p<2$. We call this positive function the real bubble. Two 
blowup sequences $\{(x_k^1,r_k^1)\}$, $\{(x_k^2,r_k^2)\}$ are said to be essentially same if
$$
0<d<\frac{r_k^2}{r_k^1}<d',\s
and\s \frac{|x_k^2-x_k^1|}{r_k^1+r_k^2}< d''.
$$
Thus, after passing to a subsequence, $\{(x_k^1,r_k^1)\}$ and
$\{(x_k^2,r_k^2)\}$ are said to be essentially different if 
$$
\frac{|x_k^2-x_k^1|}{r_k^1+r_k^2}\rightarrow+\infty,
$$
or
$$
\frac{r_k^1}{r_k^2}+\frac{r_k^2}{r_k^1}\rightarrow+\infty.
$$
That is to say, we may assume
\begin{equation}\label{two.bubble1}
B_{Rr_k^1}(x_k^1)\cap B_{Rr_k^2}(x_k^2)=\emptyset
\end{equation}
for any fixed $R$, when $k$ is sufficiently large (See Figure 3), or
\begin{equation}\label{two.bubble2}
\frac{r_k^1}{r_k^2}\rightarrow 0, \s and\s \frac{x_k^1-x_k^2}{r_k^2}\rightarrow x_{12}
\end{equation}
(See Figure 4).

\begin{figure}[htbp]
\centering
\begin{tikzpicture}[line cap=round,line join=round,>=triangle 45,x=0.7cm,y=0.7cm]
%\clip(-1.4,-0.9) rectangle (26.22256298253256,13.891445635803642);
\begin{scope}[scale=.9,xshift=1cm]
\draw [line width=0.5pt] (3.2,3.9884568426640294) circle (1.7531046320955417cm);
\draw [line width=0.5pt] (2.1190075322844923,4.354546003679194) circle (0.4401749873244177cm);
\draw [line width=0.5pt] (4.198393966850634,3.344139919277339) circle (0.5857152002936413cm);
\begin{scriptsize}
\draw (2.1190075322844923,4.354546003679194) circle (0.5pt);
\draw(2.2837459539383103,4.804223052312411) node {$(x_k^1,r_k^1)$};
\draw (4.198393966850634,3.344139919277339) circle (0.5pt);
\draw (4.35606145790617,3.779975159546688) node {$(x_k^2,r_k^2)$};
\end{scriptsize}
\end{scope}
\draw [->,line width=0.5pt] (6.863523059041039,3.9591697097828162) -- (8.459671801067163,3.9884568426640294);
\draw [shift={(10.251156797088294,2.4184685072801995)},line width=0.5pt]  plot[domain=-1.4715094205907358:1.3855409008551953,variable=\t]({1.*0.6885117459838842*cos(\t r)+0.*0.6885117459838842*sin(\t r)},{0.*0.6885117459838842*cos(\t r)+1.*0.6885117459838842*sin(\t r)});
\draw [shift={(10.892563795868051,4.55671357256839)},line width=0.5pt]  plot[domain=-0.6711870849323169:4.373853321998752,variable=\t]({1.*1.5494584557148467*cos(\t r)+0.*1.5494584557148467*sin(\t r)},{0.*1.5494584557148467*cos(\t r)+1.*1.5494584557148467*sin(\t r)});
\draw [shift={(12.889669520981704,3.1420611552556683)},line width=0.5pt]  plot[domain=-3.6637752553454193:0.09406123462798456,variable=\t]({1.*0.904257733322195*cos(\t r)+0.*0.904257733322195*sin(\t r)},{0.*0.904257733322195*cos(\t r)+1.*0.904257733322195*sin(\t r)});
\draw [shift={(14.790544884767508,3.7053043453256764)},line width=0.5pt]  plot[domain=-1.6147965047524444:3.587501358923534,variable=\t]({1.*1.1090597198179426*cos(\t r)+0.*1.1090597198179426*sin(\t r)},{0.*1.1090597198179426*cos(\t r)+1.*1.1090597198179426*sin(\t r)});
\draw [shift={(14.740615236308146,1.8930274872810233)},line width=0.5pt]  plot[domain=1.569168352691191:5.529752890436214,variable=\t]({1.*0.7042914768169102*cos(\t r)+0.*0.7042914768169102*sin(\t r)},{0.*0.7042914768169102*cos(\t r)+1.*0.7042914768169102*sin(\t r)});
\draw [->,line width=0.5pt] (16.630781874925667,3.541828066225525) -- (18.226930616951787,3.5711151991067385);

\begin{scope}[yshift=0.5cm]
\draw [rotate around={-65.47227951974193:(19.954871456943362,3.329496352836744)},line width=0.5pt] (19.954871456943362,3.329496352836744) ellipse (1.405779295425666cm and 0.79334345164333cm);
\draw [shift={(22.05417744880004,3.7195428192418576)},line width=0.5pt]  plot[domain=-0.9136621988862466:2.8934087008622957,variable=\t]({1.*1.213956795421983*cos(\t r)+0.*1.213956795421983*sin(\t r)},{0.*1.213956795421983*cos(\t r)+1.*1.213956795421983*sin(\t r)});
\draw [shift={(19.274815111127914,6.160720338936132)},line width=0.5pt]  plot[domain=4.988889299566216:5.514914054255934,variable=\t]({1.*4.896181892410446*cos(\t r)+0.*4.896181892410446*sin(\t r)},{0.*4.896181892410446*cos(\t r)+1.*4.896181892410446*sin(\t r)});
\draw [shift={(21.12827418231239,-3.69415572791606)},line width=0.5pt]  plot[domain=1.1386945654040124:2.0162312239935165,variable=\t]({1.*5.170564074989353*cos(\t r)+0.*5.170564074989353*sin(\t r)},{0.*5.170564074989353*cos(\t r)+1.*5.170564074989353*sin(\t r)});

\end{scope}
\begin{scriptsize}

\draw(10.977943183228758,5.042420236676532) node {$u_k^1=r_k^1u_k(r_k^1x+x_k^1)$};

\draw(15.14639390960089,4.113451217656459) node {$u_k^2=r_k^2u_k(r_k^2x+x_k^1)$};

\draw (20.11280520359283,3.9586230478197795) node {$u^1$};

\draw(22.30421929974275,3.8633441740741308) node {$u^2$};
\end{scriptsize}
\end{tikzpicture}
\caption{\label{fig: } }
\end{figure}

\begin{figure}[htbp]
\centering
\begin{tikzpicture}[line cap=round,line join=round,>=triangle 45,x=0.7cm,y=0.7cm]

\begin{scope}[scale=0.7,xshift=-4.5cm,yshift=1cm]
\draw [line width=0.5pt] (5.3336976289623514,6.878511996827276) circle (2.3928175575158175cm);
\draw [line width=0.5pt] (3.877874847738024,7.456852553752008) circle (0.967978966460062cm);
\draw [line width=0.5pt] (3.877874847738024,7.456852553752008) circle (0.6445806902203528cm);
\begin{scriptsize}
\draw (5.3336976289623514,6.878511996827276) circle (0.5pt);
\draw (5.951923741537067,6.908426163564762) node {$(x_k^2,r_k^2)$};
\draw(3.877874847738024,7.456852553752008) circle (0.5pt);
\draw (4.01747429251296,7.825793943514338) node {$(x_k^1,r_k^1)$};
\end{scriptsize}
\end{scope}

\draw [->,line width=0.5pt] (
1.882081527435483,5.8511996827276) -- (3.57864819693493,5.8511996827276);

\begin{scope}[scale=0.85,xshift=-3.5cm]
\draw [shift={(13.735429771046554,9.817869586727742)},line width=0.5pt]  plot[domain=-1.0217871675346695:4.018180913663478,variable=\t]({1.*0.9442446673479159*cos(\t r)+0.*0.9442446673479159*sin(\t r)},{0.*0.9442446673479159*cos(\t r)+1.*0.9442446673479159*sin(\t r)});
\draw [shift={(12.539600548068565,8.549769259555779)},line width=0.5pt]  plot[domain=-0.9110146975009004:0.7419271864749055,variable=\t]({1.*0.8026982228459426*cos(\t r)+0.*0.8026982228459426*sin(\t r)},{0.*0.8026982228459426*cos(\t r)+1.*0.8026982228459426*sin(\t r)});
\draw [shift={(13.227733133711633,6.20144818187332)},line width=0.5pt]  plot[domain=1.6847193230771775:4.528319323035958,variable=\t]({1.*1.725271834878398*cos(\t r)+0.*1.725271834878398*sin(\t r)},{0.*1.725271834878398*cos(\t r)+1.*1.725271834878398*sin(\t r)});
\draw [shift={(12.869974527270388,3.9687269850795)},line width=0.5pt]  plot[domain=-1.8291372235457581:1.4927236798315864,variable=\t]({1.*0.5382339765912753*cos(\t r)+0.*0.5382339765912753*sin(\t r)},{0.*0.5382339765912753*cos(\t r)+1.*0.5382339765912753*sin(\t r)});
\draw [shift={(14.738379271597589,8.576971908255656)},line width=0.5pt]  plot[domain=2.435116337837783:4.545610976509619,variable=\t]({1.*0.6707421765342416*cos(\t r)+0.*0.6707421765342416*sin(\t r)},{0.*0.6707421765342416*cos(\t r)+1.*0.6707421765342416*sin(\t r)});
\draw [shift={(14.169476791395411,6.167861332694797)},line width=0.5pt]  plot[domain=-1.3488231165192728:1.3147355085292793,variable=\t]({1.*1.806578298870772*cos(\t r)+0.*1.806578298870772*sin(\t r)},{0.*1.806578298870772*cos(\t r)+1.*1.806578298870772*sin(\t r)});
\draw [shift={(14.872917907875506,3.901404547845744)},line width=0.5pt]  plot[domain=2.1158583192462803:5.640639612108339,variable=\t]({1.*0.5896454891796472*cos(\t r)+0.*0.5896454891796472*sin(\t r)},{0.*0.5896454891796472*cos(\t r)+1.*0.5896454891796472*sin(\t r)});

\begin{scriptsize}

\draw(13.939006260445879,6.639198662927387) node {$u_k^2=r_k^2u_k(r_k^2x+x_k^2)$};
\draw(13.919063482620887,10.228898671425725) node {$u_k^1=r_k^1u_k(r_k^1x+x_k^1)$};
\end{scriptsize}
\end{scope}

\draw [->,line width=0.5pt] (11.07827710263165,5.8511996827276) -- (12.774557661256053,5.8511996827276);

\begin{scope}[scale=0.8, xshift=-1.5cm,yshift=.4cm]
\draw [rotate around={-79.99202019855808:(21.347748222429946,8.872789779326354)},line width=0.5pt] (21.347748222429946,8.872789779326354) ellipse (0.7127677356629385cm and 0.6707932544065303cm);
\draw [rotate around={86.56636963755203:(21.377662389167437,6.539484773802427)},line width=0.5pt] (21.377662389167437,6.539484773802427) ellipse (0.9019580768877192cm and 0.831438428999241cm);
\draw [shift={(21.546351993322027,2.3874155529975996)},line width=0.5pt]  plot[domain=0.7470202307254491:2.403962025239333,variable=\t]({1.*2.882254595913317*cos(\t r)+0.*2.882254595913317*sin(\t r)},{0.*2.882254595913317*cos(\t r)+1.*2.882254595913317*sin(\t r)});
\begin{scriptsize}
\draw(21.61697572306733,8.750844225376037) node {$u^1$};
\draw (21.736632390017274,6.798740885527313) node {$u^2$};
\end{scriptsize}

\end{scope}
\end{tikzpicture}
\caption{\label{fig: } }
\end{figure}

The goal of this section is to illustrate that the volume and the diameter of the neck region will vanish as $k\to \infty$. As a preparation, we establish the following lemmas. Set $g_{q}=dt^2+g_{\mathbb{S}^{3}}$ and $dV_{g_q}=dV_q$. 
\begin{lem}\label{diamcontrol}
Let $Q=[0,3]\times \mathbb{S}^{3}$, and $Q'=[1,2]\times \mathbb{S}^{3}$. Assume $g$ is a Riemannian metric over $Q$ and $u^2g$ is a conformal metric in the conformal class $[g]$. There exist $\varepsilon_3>0$ and $\tau_3>0$ such that if
$$
||g-g_q||_{C^2(Q)} <\tau_3 \s and \s \int_{Q}|R|^2u^4dV_g<\varepsilon_3,
$$  
then we have 
$$
\mathrm{diam}(Q',u^2g)^4\leq C(\tau_3,\varepsilon_3)\mathrm{Vol}(Q',u^2g),
$$
where $d_u$ is the distance function of $u^2g$.
\end{lem}
\proof
Assume the result is not true. There exists a sequence of conformal metrics $g_k=u_k^2h_k$  with 
$$
||h_k-g_q||_{C^2} \to 0 \s and \s \int_Q|R_k|^2u_k^4dV_{h_k}\to 0,
$$
and 
$$
\mathrm{diam}(Q',g_k)^4>k\mathrm{Vol}(Q',g_k),
$$
where $d_k$ is the distance function of $g_k$.
Without of loss generality, assume $\mathrm{Vol}(Q',g_k)=1$. If not, we can take $v_k=c_ku_k$, such that $\mathrm{Vol}(Q',c_k^2g_k)=1$. It does not influence our result. Hence $\mathrm{diam}(Q',g_k)^4\to\infty$.
By $\varepsilon$-regularity, $u_k$ converges to $v$ weakly in $W^{2,p}(Q,g_q)$ and $v$ is a positive solution  of
$$
-\Delta_{g_{q}}v+R_{g_{q}}v=0
$$
on $Q$. In the above equation, $R_{g_q}$ is the scalar curvature of $g_q$ which is a constant. By the Harnack's inequality 
and 
$$
\int_{Q'}v^4dV_q= 1,
$$
$v$ is bounded in $Q'$. Thus
$$
\mathrm{diam}(Q',v^2g_{q})\leq |v|_{L^\infty(Q')}\mathrm{diam}(Q',g_{q})<\infty.
$$
By Lemma \ref{lem4}, $\mathrm{diam}(Q',g_k)\rightarrow 
\mathrm{diam} (Q',v^2g_q)$.
We get a contradiction and complete the proof.
\endproof

We have the similar consequence over balls in $(\R^4, g_{eucl})$ as follows and omit the proof:
\begin{lem}\label{diamcontrol2}
Let $g$ be a Riemannian metric over $B_3\subset (\mathbb{R}^4, g_{eucl})$. Assume $u^2g$ is a conformal metric in the conformal class $[g]$. There exist $\varepsilon_4>0$ and $\tau_4>0$ such that if
$$
||g-g_{eucl}||_{C^2(B_3)} <\tau_4 \s and \s \int_{B_3}|R|^2u^4dV_g<\varepsilon_4,
$$  
then we have 
$$
\mathrm{diam}(B_2,u^2g)^4\leq C(\tau_4,\varepsilon_4)\mathrm{Vol}(B_2,u^2g),
$$
where $d_u$ is the distance function of $u^2g$.
\end{lem}

Now, we close this section by proving the following proposition. Let $\{(x_k,r_k)\}$ be a blowup sequence and $r_0>0$ be fixed.
\begin{pro}\label{neckvolume} 
There exists $\varepsilon_5$, such that if 
$$\mathrm{Vol}(B_{r_0}\setminus B_\frac{r_k}{r_0}(x_k),g_k)\leq 1 \s and \s \int_{B_{2t}\setminus B_t(x_k)}R_k^2dV_{g_k}<\varepsilon_5, \s
\forall t\in [\frac{r_k}{r_0},r_0],$$  
then 
$$
\lim_{r\to0}\lim_{k\to \infty}(\mathrm{diam}(B_r(x_k)\setminus B_{\frac{r_k}{r}}(x_k),g_k)+\mathrm{Vol}(B_r(x_k)\setminus B_{\frac{r_k}{r}}(x_k),g_k))=0.
$$
\end{pro}
\proof

Choose a normal chart around $0$, and set 
$$
\phi_k(r,\theta)=x_k+(-\log r,\theta).
$$
On this polar coordinate, we set $h_k(t,\theta)=\phi_k^*(g_0)(t,\theta)$ and
$g_k=v_k^2h_k$, where $v_k(t,\theta)=u_k(x_k+(e^{-t},\theta))e^{-t}$. We may assume 
$$
||h_k-g_q||_{C^2((-\log r_0,-\log (r_k/r_0))\times \mathbb{S}^3)}<\tau_1 \s(defined\s in\s Lemma \s \ref{TCL}),
$$

Since $\mathrm{Vol}((-\log r_0,\log\frac{r_k}{r_0})\times \mathbb{S}^3, g_k)\leq 1$, for any $\varepsilon>0$,
we can choose $m<\frac{1}{2\varepsilon}$, and $T$, $T_k$ as follows:

1) $-\log r_0<T<-\log r_0+m$, $-\log r_k/r_0-m<T_k<-\log r_k/r_0$.

2) $(T_k-T)/L$ is an integer, $L\in (L_0,2L_0)$, where $L_0$ is mentioned in Lemma \ref{TCL}. 

3) $\mathrm{Vol}([T,T+L]\times \mathbb{S}^3, g_k)+\mathrm{Vol}([T_k-L,T_k]\times \mathbb{S}^3, g_k)<\varepsilon.$

Set 
$$
A_k=(\int_{[T,T+L]\times \mathbb{S}^{3}}v_k^2 dV_0)^{\frac{1}{2}}+(\int_{[T_k-L,T_k]\times \mathbb{S}^{3}}v_k^2 dV_0)^{\frac{1}{2}},
$$
and
$$
A_k'=\mathrm{Vol}([T,T+L]\times \mathbb{S}^{3},g_k)^\frac{1}{4}+
\mathrm{Vol}([T_k-L,T_k]\times \mathbb{S}^{3},g_k)^\frac{1}{4}.
$$
Obviously,
$$
A_k\leq CA_k' \s and \s A_k'<\sqrt[4]{\varepsilon}.
$$
 
Firstly, we prove 
$$
\mathrm{diam}([T,T_k]\times \mathbb{S}^3,g_k)\leq C\sqrt[4]{A_k'},\s \mathrm{Vol}([T,T_k]\times \mathbb{S}^3,g_k)\leq CA_k'.
$$
$v_k$ satisifes the following equation:
$$
-\Delta_{h_k}v_k+R(h_k)v_k=R(g_k)v_k^3.
$$
Since for any $t$,
\begin{align*}
\int_{[t,t+\log 2]\times \mathbb{S}^3} R_k^2v_k^4dV_{h_k}=\int_{B_{e^{-t}}(x_k)\setminus B_{\frac{1}{2}e^{-t}}(x_k)} R_k^2u_k^4dV_{g_0}< \varepsilon_5, 
\end{align*}
we can choose $\varepsilon_5$ such that $3L\varepsilon_5/\log 2<\min\{\varepsilon_0,\varepsilon_1,\varepsilon_3\}$, then whenever $t\in [-\log r_0,-\log (r_k/r_0)-3L]$:
$$
\int_{[t,t+3L]\times \mathbb{S}^3} R_k^2v_k^4dV_{h_k}<\min\{\varepsilon_0,\varepsilon_1,\varepsilon_3\}.
$$

Putting $Q_i=[T+iL,T+(i+1)L]$, by Lemma \ref{regularity} and Lemma \ref{diamcontrol}, we have
$$
\|v_k\|_{L^6(Q_i)}<C,
$$
and
$$
\mathrm{diam}(Q_i,g_k)\leq\sqrt[4]{ \mathrm{Vol}(Q_i,g_k)}.
$$
By Three circles Lemma \ref{TCL}, there exists $\delta>0$, such that
$$
\int_{Q_i}v_k^4dx<\|v_k\|^3_{L^6(Q_i)}\|v_k\|_{L^2(Q_i)}<C(e^{-i\delta L}+e^{-(n_k-i)\delta L})A_k(L).
$$
Hence
\begin{align*}
\mathrm{Vol}([T,T_k]\times \mathbb{S}^3,g_k)\leq C A_k(L) \sum_i(e^{-i\delta L}+e^{-(n_k-i)\delta L})\leq C A_k'(L),
\end{align*}
and
\begin{eqnarray*}
\mathrm{diam}([T,T_k]\times \mathbb{S}^3,g_k)&\leq&\sum_i\mathrm{diam}(Q_i,g_k)\leq C\sum_i\sqrt[4]{\mathrm{Vol}(Q_i,g_k)}\\
&\leq& C(A_k')^\frac{1}{4}(\sum_i(e^{-i\delta\frac{L}{4}}+e^{-(n_k-i)\delta\frac{L}{4}})\leq
C\sqrt[4]{A_k'}.
\end{eqnarray*}

Then when $r$ is sufficiently small and $k$ is sufficiently large, we have
$$
\mathrm{Vol}(B_r\setminus B_\frac{r_k}{r}(x_k),g_k)\leq \mathrm{Vol}([T,T_k]\times \mathbb{S}^3,g_k)\leq C\sqrt[4]{\varepsilon}.
$$
Letting $\varepsilon\rightarrow 0$, we get
$$
\lim_{k\rightarrow+\infty}\lim_{r\rightarrow 0}\mathrm{Vol}(B_r\setminus B_\frac{r_k}{r}(x_k),g_k)=0.
$$
In the same way, we derive
$$
\lim_{k\rightarrow+\infty}\lim_{r\rightarrow 0}\mathrm{diam}(B_r\setminus B_\frac{r_k}{r}(x_k),g_k)=0.
$$
\endproof

\section{The Proof of the Main Theorem}

The goal of the final section is to finish proving Theorem \ref{main}. We take $\varepsilon=\min \{\varepsilon_i, i=0,\cdots, 5\}$, where $\varepsilon_i$ come from Lemma \ref{regularity}, Lemma \ref{TCL}, Lemma \ref{lem4}, Lemma \ref{diamcontrol} Lemma \ref{diamcontrol2} and Proposition \ref{neckvolume}, and define a finite set 
$$
\mathcal{S}_0=\{x\in M: \lim_{r\to 0}\liminf_{k\to \infty}\int_{B_r(x)} R_k^2dV_{g_k}>\varepsilon\}.
$$
The points in $\mathcal{S}_0$ are said to be bubble points. Obviously, $\mathcal{S}$ and $\mathcal{S}'$ which are defined in section 2 and section 3 are two subsets of $\mathcal{S}_0$, hence $\{u_k\}$ weakly converges to $u$ in $W^{2,p}_{loc}(M\setminus \mathcal{S}_0,g_0)$ for some $p<2$ and $\{d_k\}$ converges to $d_u$ in $C^0_{loc}(M\setminus \mathcal{S}_0, g_0)$. 
Before proving our main theorem, some preparations need to be done. At first we are concerned with the situation in which there is no real bubble at bubble points.

\begin{lem}\label{nobubble}
Let $P\in\mathcal{S}_0$. If there is no real bubble at $P$,
then 
$$
\lim_{r\rightarrow 0}\lim_{k\rightarrow+\infty}\mathrm{diam}(B_r(P),g_k)=0,\s and\s
\lim_{r\rightarrow 0}\lim_{k\rightarrow+\infty}\mathrm{Vol}(B_r(P),g_k)=0.
$$
\end{lem}
\proof
Prove by the method of induction. Let
$$
m=[\frac{2}{\varepsilon}\lim_{r\to 0}\liminf_{k\to \infty}\int_{B_r(p)} R_k^2dV_{g_k}].
$$

If $m=1$, there is no concentration phenomena at $P$, $\{u_k\}$ converges to $u$ weakly in $W^{2,p}(B_1(P))$. Then we have
\begin{align*}
\lim_{r\rightarrow 0}\lim_{k\rightarrow+\infty}\mathrm{Vol}(B_r(P),g_k)&=\lim_{r\rightarrow 0}\mathrm{Vol}(B_r(P),g_u)
\\ &=\lim_{r\rightarrow 0}\int_{B_r(P)} u^4dV_{g_0}=0
\end{align*}
By Lemma \ref{diamcontrol2}, $\lim_{r\rightarrow 0}\lim_{k\rightarrow+\infty}\mathrm{diam}(B_r(p),g_k)=0$. The proof is completed.

Assume the statements holds for $m=m_0-1$. Let's consider about $m=m_0$. Choose a normal chart around $P$, set $P=0$ and $(g_0)_{ij}=\delta_{ij}+O\{|x|^2\}$. 
Then define
$$
r_k=\inf\{r|\int_{B_r(x)}|R_k|^2dV_{g_k}=\frac{\varepsilon}{2}, \exists x\in B_\delta(0)\}.
$$
Take $x_k\in B_\delta(0)$ such that $\int_{B_{r_k}(x_k)}|R_k|^2dV_{g_k}=\frac{\varepsilon}{2}$. 
Then define
$$
t_k=\inf\{r|\int_{B_r(x_k)}|R_k|^2dV_{g_k}=(m_0-1)\frac{\varepsilon}{2}\}.
$$ 
Since $r_ku_k(x_k+r_kx)$ converges to 0 weakly in $W^{2,p}(\R^4)$, we may assume
$t_k/r_k\rightarrow+\infty$. 

Obviously, there exists a $r$, such that $\forall t\in[t_k/r,r]$,
$$
\int_{B_{2t}\setminus B_t(x_k)} R_k^2dV_{g_k}<\varepsilon/2.
$$
Then, by Proposition \ref{neckvolume}, we get
$$
\lim_{r\rightarrow 0}\lim_{k\rightarrow+\infty}
diam_{g_k}(B_r\setminus B_{\frac{t_k}{r}}(x_k),g_k)+\mathrm{Vol}(B_r\setminus B_{\frac{t_k}{r}}(x_k),g_k)=0.
$$

Let $v_k=t_ku_k(x_k+t_kx)$ and $h_k(x)=(g_0)_{ij}(x_k+t_kx)dx^i\otimes dx^j$. Obviously $h_k$ converges to the canonical Euclidean metric $g_{eucl}$ smoothly. We have
$$
-6\Delta v_k+t_k^2 R(h_k)v_k=\hat{R}_kv_k^3,
$$
where $\hat{R}_k$ denotes the rescaling scalar curvature responding to the new conformal factor $v_k$. 
Then 
$$
\lim_{k\to \infty}\int_{B_\frac{r_k}{t_k}(0)}\hat{R}_k^2v_k^4dx=\int_{B_{r_k}(x_k)}R_k^2u_k^4dV_0=\frac{\varepsilon}{2},
$$
and
$$
\lim_{k\to \infty}\int_{B_1(0)}\hat{R}_k^2v_k^4dx=\int_{B_{t_k}(x_k)}R_k^2u_k^4dV_0= (m_0-1)\frac{\varepsilon}{2}.
$$
Hence
$$
\lim_{k\to \infty}\int_{B_1\setminus B_\frac{r_k}{t_k}(0)}\hat{R}_k^2v_k^4dx= (m_0-2)\frac{\varepsilon}{2}.
$$
Moreover, for any sufficiently large $R$,
$$
\lim_{k\to \infty}\int_{B_R\setminus B_1(0)} \hat{R}_k^2v_k^4dx\leq \frac{\varepsilon}{2}
$$
Letting $R\to \infty$, we get $\forall x\in \R^4\setminus \{0\}$
\begin{align*}
\lim_{r\to0}\liminf_{k\to\infty}\int_{B_r(x)} \hat{R}_k^2v_k^4dx &\leq \lim_{r\to 0}\liminf_{k\to\infty}\int_{B_r(x)\cap( B_1(0)\setminus\{0\})} \hat{R}_k^2v_k^4dx
\\ &+ \lim_{r\to 0}\liminf_{k\to\infty}\int_{B_r(x)\cap(\R^4\setminus B_1(0))} \hat{R}_k^2v_k^4dx
\\ &\leq (m-1)\frac{\varepsilon}{2},
\end{align*}
and 
$$
\lim_{r\to0}\liminf_{k\to\infty}\int_{B_r(0)} \hat{R}_k^2v_k^4dV_0\leq \int_{B_1(0)}\hat{R}_k^2v_k^4dx=(m_0-1)\frac{\varepsilon}{2}.
$$
Denote $\mathcal{S}'(v_k)$ as a concentration set of $v_k$ over $\R^4$, i.e., 
$$\mathcal{S}'(v_k)=\{x\in \R^4: \lim_{r\to 0}\liminf_{k\to \infty}\int_{B_r(x)}\hat{R}_k^2v_k^4dx>\varepsilon\}.$$
Hence $\forall x\in \mathcal{S}'\subset \R^4$,
$$
\lim_{r\to0}\liminf_{k\to\infty}\int_{B_r(x)} \hat{R}_k^2v_k^4dV_0\leq (m_0-1)\frac{\varepsilon}{2}.
$$
Using the induction hypothesis, we obtain
$$
\lim_{r\rightarrow 0}\lim_{k\rightarrow+\infty}
diam_{\hat{g}_k}(B_r(x))+\mathrm{Vol}(B_r(x),\hat{g}_k)=0, \s \forall x \in\mathcal{S}',
$$
where $\hat{g}_k=v_k^2g_{eucl}$. The proof is finished.
\endproof

Then we present $\{(M,g_k)\}$ can not be pinched, as $k\to \infty$. 
\begin{lem}\label{pinch}
If there exists $r_k\rightarrow 0$, such that 
$$
\lim_{R\to +\infty}\liminf_{k\to \infty} \mathrm{Vol}(B_{Rr_k}(x_k),g_k)>0, \s \s \lim_{R\to +\infty}\liminf_{k\to \infty} \mathrm{Vol}(M\setminus B_{Rr_k}(x_k),g_k)>0,
$$
then the first eigenvalue: $\lambda_1(\Delta_{g_k})\to 0$, we call $(M,g_k)$ can be pinched.
\end{lem}
\proof
Let 
$$
V_1=\lim_{R\to +\infty}\liminf_{k\to \infty} \mathrm{Vol}(B_{Rr}(x_k),g_k)>0,\s V_2=1-V_1.
$$

Given an $\tau>0$, 
we assume, for sufficiently large $k$,
$$
0<V_1-\mathrm{Vol}(B_{2^{m_0}r_k}(x_k),g_k)<\tau.
$$
Since
$$
\sum_{i=m_0+1}^m \mathrm{Vol}(B_{2^{i+1}r_k}\setminus B_{2^{i-1}r_k}(x_k),g_k)\leq 2\mathrm{Vol}(B_{2^{m+1}r_k}\setminus B_{2^{m_0}r_k}(x_k),g_k)\leq 2,
$$
we can find $i_k$ between
$m_0+1$ and $m_0+\frac{1}{2\tau}$, such that
$$
\mathrm{Vol}(B_{2^{i_k+1}r_k}\setminus B_{2^{i_k-1}r_k}(x_k),g_k)<\tau.
$$
Let $t_k=2^{i_k}r_k$ and  $A_k=B_{2t_k}\setminus B_{t_k/2}(x_k)$. We have
$$
|\mathrm{Vol}(B_{t_k/2}(x_k),g_k)-V_1|<2\tau,\s
|\mathrm{Vol}(B_{2t_k}(x_k)^c,g_k)-V_2|<2\tau,\s 
\mathrm{Vol}(A_k,g_k)<\tau.
$$
Let $\eta:\mathbb{R}\to \mathbb{R}$ be a smooth decreasing function, 
which is $1/V_1$ on $(-\infty,1/2]$ and $-1/V_2$ on $[2,+\infty)$,
and set
\begin{eqnarray*}
\Psi_k(x)=
\begin{cases}
1/V_1 &x\in B_{\frac{t_k}{2}}(x_k) \cr 
\eta(\frac{d_0(x,x_k)}{t_k})&x\in A_k  \cr 
-1/V_2, & x\notin  B_{2t_k}(x_k). 
\end{cases}
\end{eqnarray*}
Then
\begin{align*}
\int_M|\nabla_k\Psi_k|^2dV_{g_k}&=\int_{A_k} |\nabla_k\Psi_k|^2dV_{g_k}
\\ 
&=\int_{A_k} |\eta'|^2\frac{|\nabla_0 d_0(x,x_k)|^2}{t_k^2}u_k^2dV_0
\\
&\leq C(\int_{A_k} \frac{|\nabla_0 d_0(x,x_k)|^4}{t_k^4}dV_0)^\frac{1}{2}\mathrm{Vol}(A_k,g_k)^{\frac{1}{2}}
\\
& =C(\int_{A_k} \frac{1}{t_k^4}dV_0)^\frac{1}{2}\mathrm{Vol}(A_k,g_k)^{\frac{1}{2}}
\\ 
&\leq C\sqrt{\tau}.
\end{align*}
Since
$$
\int_{A_k}|\Psi_k|dV_{g_k}\leq \max\{1/V_1,1/V_2\}\tau<C\tau,
$$
and
\begin{eqnarray*}
\left|\int_{M\setminus A_k}\Psi_k dV_{g_k}\right|&=&\left|\frac{\mathrm{Vol}(B_{t_k/2}(x_k),g_k)}{V_1}-
\frac{\mathrm{Vol}(M\setminus B_{2t_k}(x_k),g_k)}{V_2}\right|\\
&=&\left|\frac{\mathrm{Vol}(B_{t_k/2}(x_k),g_k))-V_1}{V_1}-
\frac{\mathrm{Vol}(M\setminus B_{2t_k}(x_k),g_k)-V_2}{V_2}\right|\\
&\leq& C\tau,
\end{eqnarray*}
we get
$$
\left|\int_M\Psi_k dV_{g_k}\right|\leq C\tau.
$$
Also,
\begin{eqnarray*}
\|\Psi_k-\bar\Psi\|_{L^2(M,g_k)}&\geq&\|\Psi_k\|_{L^2(M\setminus A_k,g_k)}-C\tau
\\
&\geq& \left(\frac{\mathrm{Vol}(B_{\frac{t_k}{2}}(x_k),g_k)}{V_1^2}+ 
\frac{\mathrm{Vol}(M\setminus B_{2t_k}(x_k),g_k)}{V_2^2}\right)^\frac{1}{2}
-C\tau
\\
&\geq& \left(\frac{V_1-2\tau}{V_1^2}+\frac{V_2-2\tau}{V_2^2}\right)^\frac{1}{2}-C\tau.
\end{eqnarray*}
Hence
$$
\lambda_1(\Delta_{g_k})\leq\frac{\int_M|\nabla \Psi_k|^2dV_{g_k}}{\int_M|\Psi_k-\bar{\Psi}_k|^2dV_{g_k}}\leq C\frac{\sqrt{\tau}}
{V_1^{-1}+V_2^{-1}-C\tau}.
$$
Letting $\tau\to 0$, we get $\lambda_1(\Delta_{g_k})\to 0$.
\endproof

\textbf{ The proof of Theorem \ref{main}} 

Finally, we finish proving our main theorem by discussing the following two cases. 

\textbf{Case 1: $u\neq 0$.}

In this case,  there is no real bubble.
We can argue by contradiction. If there exists 
a blowup sequence $\{(x_k,r_k)\}$ of $\{u_k\}$ at $P\in \mathcal{S}_0$, such that
$\{r_ku_k(x_k+r_kx)\}$ converges to a nontrivial bubble, then
$$
\lim_{r\to 0}\lim_{k\to\infty}\mathrm{Vol}(B_{\frac{r_k}{r}}(x_k))>0.
$$
Since $u\neq 0$
$$\lim_{r\to 0}\lim_{k\to\infty}\mathrm{Vol}(M\setminus B_{\frac{r_k}{r}}(x_k),g_k)\geq \lim_{r\to 0}\int_{M\setminus B_{r}(x_k)} u^4>0.$$
Hence, by Lemma \ref{pinch}, we can get a contradiction. According to Proposition \ref{equal} and Lemma \ref{nobubble}, we can get our conclusion.

\textbf{Case 2: $u=0$.}

In this case, the number of bubbles is exactly one.

\begin{lem}
There is exactly one bubble of $\{u_k\}$. 
\end{lem}
\proof Since $u=0$, we have
$$
\mathrm{Vol}(M,g_k)=\lim_{r\rightarrow 0}\sum_{P\in\mathcal{S}}\lim_{k\rightarrow+\infty}\mathrm(B_r(P),g_k).
$$
Thus, by Lemma \ref{nobubble}, there is at least one real bubble.

Assume there are two blowup sequences $\{(x_k^1,r_k^1)\}$ and $\{(x_k^2,r_k^2)\}$ of $u_k$. We assume
$\{v_k=r_k^1u_k(x_k^1+r_k^1x)\}$ and $\{w_k=r_k^2u_k(x_k^2+r_k^2x)\}$ converges to $v\neq 0$ and $w\neq 0$ weakly in $W^{2,p}_{loc}(\mathbb{R}^4\setminus \mathcal{S}_1)$ and $W^{2,p}_{loc}(\mathbb{R}^4\setminus \mathcal{S}_2)$ respectively.

If $x_k^i\rightarrow p^i$ with $p^1\neq p^2$, then
$B_{Rr_k^1}(x_k^1)\cap B_{Rr_k^2}(x_k^2)=\emptyset$ for sufficiently large $k$. Thus, 
\begin{equation}\label{pinch1}
\lim_{R\rightarrow+\infty}\lim_{k\rightarrow+\infty}\mathrm{Vol}(
B_{Rr_k^1}(x_k^1),g_k)\geq\int_{\R^4}v^4dx
\end{equation}
and
\begin{equation}\label{pinch2}
\lim_{R\rightarrow+\infty}\lim_{k\rightarrow+\infty}\mathrm{Vol}(M\setminus
B_{Rr_k^1}(x_k^1),g_k)\geq\lim_{R\rightarrow+\infty}\lim_{k\rightarrow+\infty}\mathrm{Vol}(
B_{Rr_k^2}(x_k^2),g_k)\geq\int_{\R^4}w^4dx.
\end{equation} 

If $p^1=p^2$, then one of \eqref{two.bubble1} and \eqref{two.bubble2} holds, which
also implies \eqref{pinch1} and \eqref{pinch2}.
By Lemma \ref{pinch},  $\lambda_1(\Delta_{g_k})\rightarrow 0$, which is in contradiction to the constraint (1) in Theorem \ref{main}.
\endproof

Assume $\{(x_k,r_k)\}$ is the blowup sequence, and $v_k=r_ku_k(x_k+r_kx)$ converges weakly to $v$ weakly in $W^{2,p}_{loc}
(\R^4\setminus\mathcal{S}')$, where $\mathcal{S}'$ is a finite set.
Let 
$A(x_k,r,r_k)=M\setminus B_{\frac{r_k}{r}}(x_k)\cup(\cup_{p\in\mathcal{S}'} B_{rr_k}(x_k+r_kp))$ be the neck region.  Since there is no real bubble of $\{v_k\}$ at any point in $\mathcal{S}'$, by Lemma \ref{nobubble}, for any $p\in\mathcal{S}'$, we have
$$
\lim_{r\rightarrow 0}\int_{B_r(p)}v_k^4=0,
$$
which implies that 
$$
\int_{\R^4}v^4=\lim_{r\rightarrow 0}\lim_{k\rightarrow+\infty}
\mathrm{Vol}(B_\frac{r_k}{r}(x_k),g_k).
$$
By Lemma \ref{pinch} and $\lambda_1(\Delta_{g_k})>0$,
$$
\lim_{r\rightarrow 0}\lim_{k\rightarrow +\infty}\mathrm{Vol}(M\setminus B_{r_k/r}(x_k),g_k)=0.
$$
Then
$$
\int_{\R^4}v^4dx=1.
$$
Moreover, 
\begin{align*}
\liminf_{k\to\infty}\frac{1}{6}\int_M R_kdV_{g_k}&=\liminf_{k\to\infty}\frac{1}{6}\int_M R_ku_k^4dV_0
\\&=\liminf_{k\to\infty}\frac{1}{6}\int_{M}\left(6|\nabla u_k|^2+R_0|u_k|^2)\right) dV_0
\\ &\geq \liminf_{k\to\infty}\frac{1}{6}(\int_{M\setminus A(x_k,r,r_k)}(6|\nabla u_k|^2+R_0|u_k|^2)+\int_{A(x_k,r,r_k)}R_0u_k^2)
\\&\geq\int_{\R^4}|\nabla v|^2+o(r)
\\ &\geq Y(\mathbb{S}^{4})+o(r).
\end{align*}
$Y(\mathbb{S}^{4})$ is the Yamabe constant of the standard Sphere $\mathbb{S}^{4}$. We can exclude this case, by the constraint (2) in Theorem \ref{main}.


{\small\begin{thebibliography}{2}
\bibitem{Adam} R. Adams, J. Fournier, Sobolev spaces.  Second edition. Pure and Applied Mathematics (Amsterdam), 140. Elsevier/Academic Press, Amsterdam, 2003. 

\bibitem{Aldana-Carron-Tapie}
C. Aldana, G. Carron and S. Tapie,  $A_\infty$ weights and compactness of conformal metrics under $L^{n/2}$ curvature bounds. {\em arXiv:1810.05387}.

\bibitem{B} J. Br\"uning, On the compactness of isospectral potentials. {\em Comm. Partial Differential Equations.}, {\bf 9} (1984), 687-698.

\bibitem{Chang-Yang1} S.Y.A. Chang, P. Yang,  Compactness of isospectral conformal metrics on $S^3$. {\em Comment. Math. Helv.}, {\bf 64} (1989), no. 3, 363-374.

\bibitem{Chang-Yang2} S.Y.A. Chang, P. Yang,  Isospectral conformal metrics on 
3-manifolds. {\em J. Amer. Math. Soc.}, {\bf 3} (1990), no. 1, 117-145.

\bibitem{Chen-Li} J. Chen, Y. Li, Homotopy classes of harmonic maps of stratified 2-spheres and application to geometric flows. {\em Advances in Mathematics.}, {\bf 263} (2014), 357-388.

\bibitem{Dong-Li-Xu} C. Dong, Y. Li and K. Xu, Sobolev metrics and conformal metrics with $\int_{M}|R|^\frac{n}{2}dV_g$ bounds. {\em arXiv:1907.12464}.

\bibitem{Gursky}M. Gursky, Compactness of conformal metrics with integral bounds on curvature. {\em Duke Math. J.}, {\bf 72} (1993), no. 2, 339-367.

 \bibitem{Gilkey}P. Gilkey, The spectral geometry of a Riemannian manifold. {\em J. Diff. Geom.}, {\bf 10} (1975), 601-618.

\bibitem{LGZ}Y. Li, G. Wei and Z. Zhou, John-Nirenberg Radius and collapse in conformal geometry. {\em arXiv:1708.02176}.

\bibitem{Li-Zhou}Y. Li, Z. Zhou, Conformal metric sequences with integral-bounded scalar curvature. {\em arXiv:1706.03919}.

\bibitem{L-W} X. Liu, Z. Wang, On isospectral compactness in conformal class for 4-manifolds. {\em Communications in Contemporary Mathematics.}, {\bf 21} (2019), no. 5, 1850041, 26 pp.

\bibitem{Moser} J. Moser, On Harnack's theorem for elliptic differential equations. {\em Communications on Pure and Applied Mathematics.}, {\bf 14} (1961), no.3, 577-591.

\bibitem{OPS} B. Osgood, R. Phillips and P. Sarnak, Compact isospectral sets of surface. {\em J. Funct. Anal.}, {\bf 80} (1988), 212-234

\bibitem{Parker} T. Parker, Bubble tree convergence for harmonic maps. {\em J. Differential Geometry.} {\bf 44} (1996), 595-633

\bibitem{Xu1} X. Xu, On compactness of isospectral conformal metrics of 4-manifolds. {\em Nagoya Math. J.}, {\bf 140} (1995), 77-99

\bibitem{Xu2} X. Xu, On compactness of isospectral conformal metrics of 4-sphere. {\em Comm. Anal. Geom.}, {\bf 3} (1995), 335-370.


\end{thebibliography}}
\end{document}